\documentclass{article}
\usepackage{fullpage}
\usepackage{graphicx}
\usepackage{subfigure}
\usepackage{booktabs} 
\usepackage{adjustbox,lipsum}
\usepackage{caption}
\usepackage{wrapfig}
\usepackage{multirow}
\usepackage[parfill]{parskip}

\usepackage{url}
\usepackage[colorlinks=True,linkcolor=magenta,citecolor=blue,urlcolor=blue,pagebackref=true,backref=true]{hyperref}

\usepackage{amsmath}
\usepackage{amssymb}
\usepackage{mathtools}
\usepackage{amsthm}

\usepackage[utf8]{inputenc} 
\usepackage[T1]{fontenc}    
\usepackage{url}            
\usepackage{booktabs}       
\usepackage{amsfonts}       
\usepackage{nicefrac}       
\usepackage{microtype}      
\usepackage{xcolor}         
\usepackage{enumitem}       

\usepackage{algorithm}
\usepackage{algpseudocode}
\usepackage{hoang}

\begin{document}

\title{Stochastic Constrained Decentralized Optimization for Machine Learning with Fewer Data Oracles: a Gradient Sliding Approach}

\author{ Hoang Huy Nguyen \textsuperscript{1}\thanks{\noindent\textsuperscript{1} H. Milton Stewart School of Industrial \& Systems Engineering, Georgia
Institute of Technology, Atlanta, GA, 30332, USA}, Yan Li \textsuperscript{1}, Tuo Zhao \textsuperscript{1}}
\date{}


\maketitle

\begin{abstract}
  In modern decentralized applications, ensuring communication efficiency and privacy for the users are the key challenges. In order to train machine-learning models, the algorithm has to communicate to the data center and sample data for its gradient computation, thus exposing the data and increasing the communication cost. This gives rise to the need for a decentralized optimization algorithm that is communication-efficient and minimizes the number of gradient computations. To this end, we propose the primal-dual sliding with conditional gradient sliding framework, which is communication-efficient and achieves an $\varepsilon$-approximate solution with the optimal gradient complexity of $O(1/\sqrt{\varepsilon}+\sigma^2/{\varepsilon^2})$ and $O(\log(1/\varepsilon)+\sigma^2/\varepsilon)$ for the convex and strongly convex setting respectively and an LO (Linear Optimization) complexity of $O(1/\varepsilon^2)$ for both settings given a stochastic gradient oracle with variance $\sigma^2$. Compared with the prior work \cite{wai-fw-2017}, our framework relaxes the assumption of the optimal solution being a strict interior point of the feasible set and enjoys wider applicability for large-scale training using a stochastic gradient oracle. We also demonstrate the efficiency of our algorithms with various numerical experiments.
\end{abstract}

\section{Introduction}
With growing demands for efficient learning algorithms for big-data applications, it is critical to consider decentralized algorithms that allow the workers to cooperate and leverage the combined computational power \cite{diffusion-strategies-learning-over-networks}. In this work, we will investigate the decentralized  optimization problem where workers will collaboratively in a decentralized manner to minimize the following constrained objective function:
\begin{equation}
    \label{eqn: problem}
    \min_{x \in \mathcal{X}} \sum_{i = 1}^m f_i(x).
\end{equation}
Here, $f_i$ are smooth, convex functions defined over a closed convex set $\mathcal{X} = \mathcal{X}_1 \times ... \times \mathcal{X}_m$ where $f_i: \mathcal{X}_i \rightarrow \mathbb{R}$ and $\mathcal{X}_i$ is a non-empty closed convex set in $\mathbb{R}^d$. Given that many large-scale machine learning applications involve training large-scale models with a large number of parameters, we also assume that the local objectives are high-dimensional functions, i.e., $d$ is large. Under this setting, each worker collects data locally and performs numerical operations using the said local data, and then passes information to the neighboring workers in a communication network. For this reason, no worker has full knowledge about other workers’ local objectives or the communication network, which helps preserve the privacy of the local data. Thus, in these decentralized and stochastic optimization problems, the workers must communicate with their neighboring workers to propagate the distributed information to every location in the network. Each network worker $i$ is associated with the local objective function $f_i(x)$ and all workers will cooperatively minimize the system objective $f(x)$ as the sum of all local objective $f_i$’s without having the full knowledge about the global problem and network structure. Decentralized optimization has many practical applications in signal processing, control, and machine learning among many other disciplines \cite{durham-pursuit-evasion, Consensus-Based-SVM, mobile-auto-agents}. 


In this work, the communication network between the workers will be described as a connected undirected graph $G = (\mathcal{N},\mathcal{E})$ where $\mathcal{N}$ is the set of indices of workers and $\mathcal{E} \subset \mathcal{N} \times \mathcal{N} = \{1,...,m\}$ is the set of edges between them where each pair of connected workers will communicate directly with each other via the edge connecting them. For convenience, we assume that there exists a loop $(i, i)$ for all workers $i \in \mathcal{N}$. The goal of the workers is to sample the first-order information of their own local objective to solve the decentralized optimization problem \ref{eqn: problem} which can be rewritten as a linearly constrained optimization problem:
\begin{wrapfigure}{r}{0.4\columnwidth}
    \vspace{-0.05in}
    \centering
    \includegraphics[width = 0.4\columnwidth]{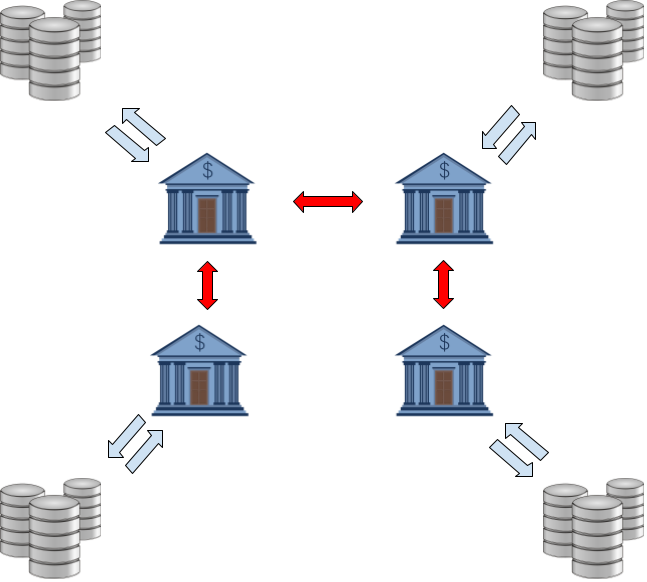}
    \caption{The optimization model where the data has to be queried and the optimizer has to communicate with local workers. Grey arrows represent local data oracle access and red arrows represent the communication between the nodes.}
    \label{fig:model_diagram}
\end{wrapfigure}
\begin{equation}
    \label{eqn: primal-dual problem} 
    \min_{x \in \cX} \sum_{i = 1}^m f_i(x) \text{ s.t } \cA x = 0.
\end{equation}

Here, we have $x = \br{x^{(1)},...,x^{(m)}} \in \cX$ with $x^{(i)} \in \cX_i \subset \R^d$ is the $d$-dimensional vector at node $i$ for $i \in \{1,...,m\}$. We also have the linear constraint matrix $\cA = \cL \otimes I_{d} \in \R^{md \times md}$. Here, the operator $\otimes$ is the Kronecker product, $I_d$ is a $d \times d$ identity matrix and $\cL$ is an $m \times m$ Laplacian matrix of the graph $G$. 
Furthermore, in order to handle the constraints in the optimization problems, one may resort to projection-based methods to solve \eqref{eqn: primal-dual problem}. However, employing such projection-based methods can be extremely computationally prohibitive \cite{complexity-LO-oracles}, as we recall that the local objectives are high-dimensional functions. To overcome this issue, projection-free methods such as the Frank-Wolfe algorithm \cite{frank-wolfe} and the conditional gradient sliding method \cite{Lan2016ConditionalGS} can be utilized to find solutions in the constraint set via a linear oracle, which is usually cheaper than using a direct projection.

To motivate our constrained decentralized optimization framework, we consider a practical setting in that the optimizer accesses the data from third-party data storage via a data oracle as described in Figure \ref{fig:model_diagram} and the workers are connected via a decentralized topology. One example of such a setting is an internal banking system of several branches, in which each local branch connects to local data storage via an internal communication network, and the branches communicate with each other via some outer communication network that is not necessarily a complete graph. Since the local worker has to communicate with the data server in order to sample the data for the computation of the gradient, this process exposes the data to the outer network (which increases the risk of a data breach) as well as causes the gradient computation to be more expensive. The aforementioned concerns motivate us to develop an efficient constrained decentralized optimization algorithm in terms of gradient sampling/data oracle and communication complexity.



\subsection{Related works}
To solve problem \ref{eqn: problem}, most algorithms rely on near-neighbor information to update their iterates as the graph $G$ is not necessarily a complete graph. Thus, developing decentralized communication-efficient algorithms or computation complexity independent of the graph topology is critical for modern big-data applications \cite{nedic2018network}. Prior works on communication efficient methods include \cite{koloskova2019decentralized,lan2017communicationefficient,chen2018lag}. Regarding communication protocols among workers, a popular approach is to use an average consensus protocol that takes a weighted average of the local updates in order to reach consensus across all workers \cite{dimakis-gossip, nedic-subgradient-projection, consensus-based-practical-issues-konstantinos, Duchi_2012_distributed, wai-fw-2017}. Some other approaches include a second-order approach such as \cite{mokhtari2016decentralized, mokhtari2016dqm}, and ADMM-based decentralized optimization \cite{Shi_2014_ADMM_decentralized, Terelius2011DecentralizedMO, aybat-distributed-linearized-admm, boyd-admm, ermin-wei-distributed-admm}. Recently, a few decentralized algorithms could match the complexities of those of centralized algorithms for the case when the feasible region is $\R^d$ \cite{kovalev-richtarik-optimal-decentralized}, \cite{huan-li-decentralized-accelerated-increasing-penalty}, \cite{huan-li-decentralized-extra-diging} but none of these could be applied to a problem with a general constraint set $\cX$ as \eqref{eqn: problem}. Lastly, only a few prior works have discussed stochastic decentralized Frank-Wolfe methods such as \cite{decen-stoch-fw-ijcai2021}, \cite{decen-tracking-submodular-pmlr-v89-xie19b}, but these works focus on tackling the decentralized DR-submodular optimization problem rather than the constrained decentralized convex optimization problem, which is the problem of our interest. Another work on stochastic decentralized Frank-Wolfe is \cite{Xian_Huang_Huang_decen_fw_2021} where the authors use stochastic sampling and quantization to obtain the gradient estimates, but this yields a sub-optimal convergence guarantee for the convex and strongly-convex case.

Note that most prior works, with the exception of \cite{wai-fw-2017}, are projection-based such that each iteration of the above algorithm requires a projection step onto the constraint set $\mathcal{C}$. Since the projection step may be computationally costly for high-dimensional problems \cite{complexity-LO-oracles}, a common approach to developing a projection-free algorithm is to adopt the framework of the Frank-Wolfe algorithm (FW) where the algorithm will solve a cheap linear optimization subroutine instead of using an expensive projection method to obtain a solution in the constraint sets \cite{frank-wolfe}. Since in some cases it is cheaper to solve the linear optimization problem: $\arg \min_{x \in \mathcal{X}} \langle g, x \rangle$ than to perform a projection onto $\mathcal{X}$ (here $g$ is the sampled gradient from a gradient oracle), the conditional gradient methods are in the interests of the optimization and machine learning community at large \cite{Beck2004ACG}. The FW algorithm has always gained a lot of interests \cite{freund2014new, pmlr-v28-jaggi13, Projection-Free-Adaptive-Gradients, frank-wolfe-strongly-convex, haipeng-vr-sfw} and has many applications, such as in semidefinite programming \cite{yurtsever2018conditional}, low-rank matrix completion (with applications in recommender systems \cite{matrix_recommend, rank-minimization-system-approx}) and robust PCA \cite{robust_PCA}. In the constrained decentralized optimization setting, the FW algorithm also has a couple of practical applications such as electric vehicle charging \cite{fw-electric-vehicle-charging} and traffic assignment \cite{fw-traffic-assignment}. However, the FW methods use a gradient call for every LO subroutine, which means that its gradient sampling complexity could be improved further. Recently, a condition gradient sliding (CGS) method is proposed in \cite{Lan2016ConditionalGS} that is able to solve an $\varepsilon$-approximate solution of constrained optimization problems with an improved gradient sampling complexity $O\left( 1/\sqrt{\varepsilon} \right)$ gradient evaluations and $O\left( 1/\varepsilon \right)$ LO oracle calls. Several extensions to the conditional gradient sliding problem include extensions to non-convex problems \cite{nonconvex-sliding}, weakly convex setting \cite{universal-cgs, nazari2020backtracking} where a backtracking line search scheme is incorporated to the conditional gradient sliding algorithm and a second-order variant of the conditional gradient sliding \cite{carderera2020secondorder}. For most machine learning applications, the vanilla CGS is sufficient and it will be focus of this work. 


\subsection{Our Contributions}
We propose a new projection-free decentralized optimization method, which enjoys gradient sampling and communication efficiency. Our proposed method leverages an inexact primal-dual sliding framework (I-PDS), which is inspired by \cite{lan2021graph} for convex decentralized optimization. Different from \cite{lan2021graph}, which assumes each constrained subproblem can be solved exactly, our I-PDS framework only requires the constrained subproblem to be solved approximately, which can be done by applying the conditional gradient sliding method in \cite{Lan2016ConditionalGS}. Compared to the prior work \cite{wai-fw-2017}, our method leads to a significant reduction in terms of data oracle calls. Our contributions can be summarized as follows:



First, we incorporate the conditional gradient sliding method to solve the linear optimization subproblem while achieving better sampling complexity than \cite{wai-fw-2017} for large data regimes, and consequently, requiring fewer accesses to the data oracle. It is noticeable that the gradient sampling complexity is invariant to the graph topology, i.e. our algorithm is independent of the spectral gap of the graph Laplacian, and in the same order as those of centralized methods for solving stochastic and deterministic problems. Note that any naive combination of the primal-dual sliding with any Frank-Wolfe-like projection-free method will \emph{not} yield the optimal gradient sampling complexity since each Frank-Wolfe step requires a gradient sample. To achieve the optimal gradient sampling complexity, we have to be able to compute multiple linear optimization steps using a single gradient sample, which is done by cleverly leveraging the gradient sliding method.

Second, as our algorithm also admits a stochastic gradient oracle while the consensus-based Decentralized Frank-Wolfe algorithm \cite{wai-fw-2017} requires an exact gradient oracle, it implies that our algorithm is more robust to noise and more versatile for machine learning applications as it allows the stochastic approximation of the gradient of $f_i$. As shown in Table \ref{complexity-table}, our gradient sampling complexity is better when $M/\rho >> \sigma^2/\varepsilon$ where $M$ is the number of training data points and $\rho$ is the spectral gap of the communication graph. 

Third, we provide theoretical analysis for the LO complexity of the decentralized communication sliding algorithm for the convex and strongly convex functions. We obtain the complexity $O \left( 1/\varepsilon^2 \right)$ for the LO complexity for both convex and strongly-convex settings as well as for both deterministic and stochastic settings. Our analysis does not require the assumption the optimal solution lies strictly in the interior of the feasible region, in contrast to \cite{wai-fw-2017}, where the authors require $\delta = \min_{x \in \overline{\cC}} \norm{x-x^*}_2 > 0$ in order to obtain the complexity $O\br{1/\sqrt{\varepsilon}}$ in the strongly-convex case. Such an assumption does not hold for many practical applications, as demonstrated in Section \ref{sec:experiments}.






\begin{table}[htb!]
\centering
\begin{adjustbox}{max width=\textwidth}
\begin{tabular}{||c c c c c c c||} 
 \hline
 \multirow{2}{*}{Algorithm} & \multirow{2}{*}{Sampling} & \multirow{2}{*}{Communication} & \multirow{2}{*}{LO} &$x^*$ strictly & Topology & Stochastic\\ 
 & & & &inside $\cX$? & invariant? & training? \\
 [0.5ex] 
 \hline
 \multirow{2}{*}{DeFW (convex)} & \multirow{2}{*}{$O\br{M\rho^{-1} \varepsilon^{-1}}$} & \multirow{2}{*}{$O\br{\varepsilon^{-1}}$} & \multirow{2}{*}{$O\br{\varepsilon^{-1}}$} & \multirow{2}{*}{\checkmark} & \multirow{2}{*}{\xmark} & \multirow{2}{*}{\xmark} \\
 & & & & & & \\
 \hline
 \multirow{2}{*}{DeFW ($\mu$-SC)} & \multirow{2}{*}{$O\br{M \rho^{-1} \varepsilon^{-0.5}}$} & \multirow{2}{*}{$O\br{\varepsilon^{-0.5}}$} & \multirow{2}{*}{$O\br{\varepsilon^{-0.5}}$} & \multirow{2}{*}{\checkmark} & \multirow{2}{*}{\xmark} & \multirow{2}{*}{\xmark} \\
 & & & & & & \\
 \hline 
 \multirow{2}{*}{I-PDS (convex)} & \multirow{2}{*}{$O\br{\varepsilon^{-1} + \sigma^2\varepsilon^{-2}}$} & \multirow{2}{*}{$O\br{\varepsilon^{-1}}$} & \multirow{2}{*}{$O\br{\varepsilon^{-2}}$} & \multirow{2}{*}{\xmark} & \multirow{2}{*}{\checkmark} & \multirow{2}{*}{\checkmark} \\
 & & & & & & \\
 \hline
 \multirow{2}{*}{I-PDS ($\mu$-SC)} & \multirow{2}{*}{$O\br{\log \varepsilon^{-1} + \sigma^2\varepsilon^{-1}}$} & \multirow{2}{*}{$O\br{\varepsilon^{-0.5}}$} & \multirow{2}{*}{$O\br{\varepsilon^{-2}}$} & \multirow{2}{*}{\xmark} & \multirow{2}{*}{\checkmark} & \multirow{2}{*}{\checkmark} \\
 & & & & & & \\
 \hline
\end{tabular}
\end{adjustbox}
\caption{ Complexity comparison of algorithms. $\varepsilon$ denotes the desired accuracy, $\sigma^2$ denotes the variance of the stochastic gradient estimate, $\rho$ denotes the spectral gap, and $M$ denotes the number of training data points. \label{complexity-table}}
\end{table}

\subsection{Organization of the paper}
The rest of the paper is organized as follows: Section \ref{sec:PDS} is dedicated to describing and explaining the intuition behind the algorithms. Section \ref{sec:main-reults} covers the convergence guarantees of the algorithm. Lastly, we will demonstrate the algorithm's empirical performances in Section \ref{sec:experiments}. 


\section{Algorithms}
\label{sec:PDS}
In this section, we present the Inexact Primal-Dual Sliding (I-PDS) optimization framework and the Conditional Gradient Sliding (CGS) algorithm to solve the constrained optimization subproblem. This approach exploits the structure of the optimization problem \eqref{eqn: primal-dual problem} and allows us to both handle the linear constraint and "slides" through multiple inner updates using a single gradient call, thus gives a gradient complexity that matches the optimal bounds \cite{arjevani2015lower} and graph topology invariance as will be proved in Section \ref{sec:main-reults}.

Before delving into the algorithm, we will establish key notations essential for its comprehension. Suppose that for each $i \in \{1,...,m\}$, the local function $f_i(x)$ can be written as 
\begin{align}
    \label{eqn: f-components}
    f_i(x) = \tilde{f}_i(x) + \mu \nu_i(x)
\end{align}
where $\mu \geq 0$ is the strongly-convex parameter and $\nu_i$ is a strongly-convex function with strong convexity modulus $1$. For notational convenience, we also denote 
\begin{align*}
f(x) = \sum_{i=1}^m f_i(x),~~\tilde{f}(x) = \sum_{i=1}^m \tilde{f}_i(x),~~\nu (x) = \sum_{i=1}^m \nu_i(x)
\end{align*}
and any variable any local variable at the node $i$ with a superscript $(i)$. We reformulate \eqref{eqn: primal-dual problem} as a saddle point problem based on the Lagrangian multiplier method:
\begin{align}
    \label{eqn: saddle point}
    \min_{x \in \cX} \max_{z \in \R^{md}} f(x) + \langle \cA x, z \rangle,
\end{align}
where $z$ is the Lagrangian multiplier. We then consider the convex conjugate of $\tilde{f}$, i.e., $\tilde{f}^*(y) = \underset{x \in \cX}{\max} \langle x,y \rangle - \tilde{f}(x)$, and further convert \eqref{eqn: saddle point} into the following problem:
\begin{equation}
    \label{eqn: saddle point dual}
    \min_{x \in \cX} \max_{y, z \in \R^{md}} \mu \nu(x) + \langle x, y + \cA^Tz \rangle - \tilde{f}^*(y).
\end{equation}
At each agent $i$, we define $\forall \hat{x}^{(i)}, x^{(i)} \in \cX_i$:
\begin{align}
    \label{eqn: divergence V}
    &V_i(\hat{x}^{(i)},x^{(i)}) := \nu (x^{(i)}) - \nu(\hat{x}^{(i)}) - \langle \nu'(\hat{x}^{(i)},x^{(i)}-\hat{x}^{(i)}) \rangle,
\end{align}
where $(\tilde{f}_i^*)'$ is the subgradient of $\tilde{f}_i^*$. For the whole network, we denote them as $V(\hat{x}, x) = \sum_{i=1}^m V_i(\hat{x}^{(i)}, x^{(i)})$. Next, we will introduce the I-PDS algorithm. As shown in Algorithm \ref{alg:PDS}, the I-PDS algorithm is a double-loop algorithm. At each outer iteration, we perform the accelerated gradient descent updates. Here, the steps \eqref{eq:txk_agent} and \eqref{eq:hxk_agent} represent the acceleration steps with $\lambda_k, \tau_k$ being the step sizes. Then, at each node of the communication network, we compute an unbiased gradient sample by uniformly sampling a batch of data at each node at step \eqref{eq:yk_agent} and then aggregating these values to obtain a gradient sample vector $v_k$ for the entire network. Note that at every outer iteration, we access the data oracle only once to compute the gradient sample, so that the inner sliding procedure does not need to access new data. In contrast to the consensus-based method in \cite{wai-fw-2017}, the computations done in the outer iteration do \emph{not} involve the communication matrix $\cA$, which explains the graph topology invariance of the gradient complexity. When $\sigma = 0$, we obtain a deterministic (exact) version of the I-PDS algorithm, which can be achieved by computing the full gradient of $\tilde{f}$. In such a case, the algorithm may achieve the $\varepsilon$-approximation solution in fewer iterations but might incur more cost per iteration given the high dimensionality of the data. 





For each inner iteration, we execute update steps called communication sliding at \eqref{eq:tukt}, \eqref{eq:zkt} to communicate the first-order information between workers. In essence, the communication sliding procedure is a type of primal-dual method when it is applied to the saddle point formulation \eqref{eqn: saddle point}. At each primal-dual step, only the computation containing the matrix $\cA$ involves communication among the workers, while the rest can be done separately at each node. This idea also allows us to save communication rounds to give us a communication-efficient method. Then, at step \eqref{eq:xkt}, we solve the constrained subproblem:
\begin{align}
    \label{linear-problem}
    x_i^t = \argmin_{x \in \mathcal{X}} \mu \nu(x) + \langle y_k + &\mathcal{A}^Tz_k^t, x\rangle + \eta_k^t V(x_k^{t-1}, x) + p_k V(x_{k-1}, x).
\end{align}


In \cite{lan2021graph}, the authors assume that this optimization subroutine has a computationally efficient exact solution. However, while this assumption may hold for simple constraints such as Euclidean balls (which has a closed-form solution), it may not be practical for more complex problems involving constraints like a nuclear norm or polytope projections (e.g., flow polytopes or Birkhoff polytopes). This necessitates the use of projection-free algorithms to obtain solutions that lie within the constraint set. 

While the Frank-Wolfe algorithm is a popular method for the projection-free approach, it uses a gradient computation for each LO oracle call, which means that the gradient complexity cannot be better than the LO complexity. Thus, we opt for a gradient sliding approach in \cite{Lan2016ConditionalGS}, which allows several linear oracle updates per gradient computation. By doing so, the CGS algorithm effectively reduces the number of gradient computations, and consequently, the number of data oracle access. The main idea behind the CGS method is instead of applying the conditional gradient procedure directly to the original convex programming problem, we apply it to the subproblems of the accelerated gradient method. By carefully choosing the accuracy threshold $\eta$ for solving these subproblems, we can show that the CGS method can achieve a gradient complexity and LO complexity that match the lower bounds for smooth, convex optimization problems \cite{arjevani2015lower}, \cite{guzman2018lower}.


In our algorithmic framework, we employ an inner procedure \ref{alg: pds-cgs} to approximate a solution to the LO subproblem, ensuring that the following stopping condition is satisfied:
\begin{align}
    \label{eqn: CGS-stopping-condition}
    S(u_t) := \max_{x\in X}\:\langle g + \beta(u_t-u), u_t - x \rangle \leq \varepsilon_k^i.
\end{align}
Here, the LHS of \eqref{eqn: CGS-stopping-condition} is equivalent to $\max_{x\in X} \langle \phi'(u_t), u_t-x \rangle$ where $\phi(x) = \langle g,x \rangle + \beta \norm{x-u}_2^2/2$. Note that the inner product $\langle \phi'(u_t), u_t-x \rangle$ is often known as the Wolfe gap, and the procedure terminates once the gap is smaller than the tolerance $\varepsilon_k^i \geq 0$ where $\varepsilon_k^i$ is the tolerance at the $i$-th inner iteration of the $k$-th outer iteration of Algorithm \ref{alg:PDS}. The inner procedure can be iteratively solved by some efficient linear optimization routine until the condition \eqref{eqn: CGS-stopping-condition} is met. In contrast to the Frank-Wolfe method, the step size $\alpha_t$ in \eqref{def:alphat} is slightly easier to compute than that of the Frank-Wolfe method, which is usually done by some line search routine. In Section \ref{sec:main-reults}, we will specify suitable parameters so that we obtain competitive rates for our algorithm.




\begin{algorithm}[htb!]
	\caption{\label{alg:PDS}The I-PDS algorithm framework}
	\footnotesize
	\begin{algorithmic}
		\State Choose $x_0, \hx_0\in X^{(i)}$, and set $\hx_0 = x_{-1} = x_0, y_0=\nabla \tf(\hx_0), z_0=0$.
		\For {$k=1,\ldots,N$}
		\State At each node $i \in \{1,...,m\}$, compute:
            \begin{align}
			\label{eq:txk_agent}
			&\tx_k^{(i)} = x_{k-1}^{(i)} + \lambda_k(\hx_{k-1}^{(i)} - x_{k-2}^{(i)})
			\\
			\label{eq:hxk_agent}
			&\hx_k^{(i)} = (\tx_k^{(i)} + \tau_k\hx_{k-1}^{(i)})/(1+\tau_k)
			\\
			\label{eq:yk_agent}
			\text{Sample $v_k^{(i)}$ s.t. } &\E[v_k^{(i)}] = \nabla \tf_i(\hx_k^{(i)}) \text{ with batch size } c_k
		\end{align}
		  \State Set $x_k^{0,(i)} = x_{k-1}^{(i)}$, $z_k^{0,(i)} = z_{k-1}^{(i)}$, and $x_k^{-1,(i)} = x_{k-1}^{T_{k-1},(i)}$ (set $x_1^{-1,(i)}=x_0^{(i)}$).
		\For {$t=1,\ldots,T_k$}
		\begin{align}
		\label{eq:tukt}
		\tu_k^t = & x_k^{t-1} + \alpha_k^t(x_k^{t-1} - x_k^{t-2})
		\\
		\label{eq:zkt}
            z_k^t =& z_k^{t-1} + \frac{\cA \tu_k^t}{q_k^t}
		\\
		\label{eq:xkt}
            \nonumber
		x_k^t \approx &\argmin_{x\in \cX} \mu\nu(x) + \langle v_k + \cA^\top z_k^t, x\rangle \\
        &+ \eta_k^tV(x_{k}^{t-1}, x) + \tfrac{p_k}{2}\|x_{k-1} - x\|_2^2
		\end{align} 
		\EndFor
		\State Set $x_k = x_k^{T_k}$, $z_k = z_k^{T_k}$, $\hx_k =\sum_{t=1}^{T_k}x_k^t/T_k$, and $\hz_k = \sum_{t=1}^{T_k}z_k^t/T_k$.
		\EndFor 
		\State Output $x_N:=\left(\sum_{k=1}^{N}\beta_k\right)^{-1}\br{\sum_{k=1}^{N}\beta_k\hx_k}$.
	\end{algorithmic}
\end{algorithm}

\begin{algorithm}[htb!]
\caption{CGS procedure at the $i$-th inner iteration of the $k$-th outer iteration of Algorithm \ref{alg:PDS}}
\label{alg: pds-cgs}
\footnotesize
\begin{algorithmic}
    \Procedure{CGS}{$g,u,\beta$}
        \State Initialize $u_0 = u$ and $t = 1$.
        \While{$S(u_{t-1}) > \varepsilon_k^i$}
        \State Compute $v_t$ such that it is the optimal solution for the subproblem
            \begin{align}
                \label{def:vt}
                S(u_t) := \max_{x\in X}\: \langle g + \beta(u_t - u), u_t-x \rangle
            \end{align}
        \State Compute:
             \begin{align}
                \label{def:alphat}
                \alpha_t &= \min \left\lbrace 1, \frac{\langle \beta(u-u_t)-g, v_t-u_t \rangle}{\beta \norm{v_t-u_t}^2} \right\rbrace \\
                \label{def:ut}
                u_t &= (1-\alpha_t)u_{t-1} + \alpha_t v_t \\
                t &\leftarrow t+1
            \end{align}
        \EndWhile
        \State Return $u_{t-1}$
    \EndProcedure
\end{algorithmic}
\end{algorithm}



\section{Main results}
\label{sec:main-reults}

Before we proceed to the main results, we will present the preliminaries of the optimization problem. Recall that our decentralized optimization problem \eqref{eqn: problem} is equivalent to the linearly-constrained optimization problem \eqref{eqn: primal-dual problem}. Our goal is to achieve an $\varepsilon$-approximation solution $\overline{x}$ such that $f(\overline{x})-f^* \leq \varepsilon$ (the primal gap is within $\varepsilon$) and $\norm{\cA x} \leq \varepsilon$ (the consensus gap is within $\varepsilon$). First, we state some definitions and assumptions:
\begin{definition}
\label{def:convexity}
A function $f: \R^d \rightarrow \R$ is $\mu$-strongly-convex if its gradient exists everywhere and there exists $\mu \geq 0$ such that $\forall~x,y \in \R^d$: 
\begin{align*}
     f(y) - f(x) - \langle \nabla f(x), y-x \rangle \geq \frac{\mu \norm{x-y}_2^2}{2}.
\end{align*}
If $\mu = 0$ then we say $f$ is convex. 
\end{definition}
\begin{definition}
\label{def:smoothness}
A function $f: \R^d \rightarrow \R$ is $L$-smooth if its gradient exists everywhere and there exists $L \geq 0$ such that $\forall~x,y \in \R^d$: 
\begin{align*}
    f(y) - f(x) - \langle \nabla f(x), y-x \rangle \leq \frac{L \norm{x-y}_2^2}{2}.
\end{align*}
\end{definition}
From these definitions, we have the following assumptions for each worker $i$:
\begin{assumption}
\label{assumption: agent-smoothness}
(Worker's smoothness and convexity assumption) For each $i \in \{1,...,m\}$, the local function $f_i(x)$ can be written as $f_i(x) = \tilde{f}_i(x) + \mu \nu(x)$ where $\mu \geq 0$ is the strongly-convex parameter, $\tilde{f}_i$ is a convex, $\tilde{L}$-smooth function and $\nu$ is an $1$-strongly-convex function with respect to some norm. 
\end{assumption}
Assumption \ref{assumption: agent-smoothness} implies that $f_i$ is $\mu$-strongly convex and has a gradient Lipschitz component. It is a general assumption for the smooth, convex decentralized optimization setting \cite{lan2020_first_order_book}. Furthermore, regarding the gradient sampling process, recall that we have the following assumptions for each worker $i$:
\begin{assumption}
\label{assumption: stochastic-gradient-assumptions}
(Stochastic gradient oracle assumptions) The first-order information of $f_i$ obtained at each worker $i$ at the $k$-th outer iteration of Algorithm \ref{alg:PDS} via a stochastic oracle $G_i(x^{(i)},\xi^{(i)})$ satisfies $\forall x^{(i)} \in X^{(i)}$:
\begin{align}
    \label{eqn: gradient-assumptions}
    \E[G_i(x^{(i)},\xi^{(i)})] &= \nabla f_i(x^{(i)}), \\
    \E\left[\norm{G_i(x^{(i)},\xi^{(i)})-\nabla f_i(x^{(i)})}_2^2\right] &\leq \frac{\sigma^2}{c_k}.
\end{align}
\end{assumption}
This noise assumption is common for the optimization literature \cite{pmlr-v80-nguyen18c} and in various learning settings such as Federated Learning \cite{li2020fedavgconvergence} and Model Agnostic Meta-Learning \cite{mishchenko2023mamlconvergence}. This assumption allows us to control the variance by adjusting the batch size $c_k$. With these assumptions and setup, we have the following convergence guarantee results:


\begin{theorem}
\label{thm:bigtheorem}
Denote $N$ as the pre-determined number of outer iterations, $\tau:=2\sqrt{\tL/\mu}$ and $\Delta:=\lceil 2\tau + 1\rceil$ if $\mu>0$, and $\Delta:=+\infty$ if $\mu=0$. Suppose that the Assumptions \ref{assumption: agent-smoothness}, \ref{assumption: stochastic-gradient-assumptions} hold, $V(\cdot,\cdot) = \norm{\cdot-\cdot}_2^2$, and that the parameters in Algorithm \ref{alg:PDS} are set to the following: 

        For all $k\ \leq \Delta$:
	\begin{align}
		\label{eq:par_beforeDelta_S}
		\begin{aligned}
			&\tau_k = \tfrac{k-1}{2},\ \lambda_k = \tfrac{k-1}{k},\ \beta_k = k,\ p_k = \tfrac{4\tL}{k}, T_k = \left\lceil \tfrac{kR\|\cA\|}{\tL}\right\rceil, c_k = \left\lceil \tfrac{\min\{N,\Delta\} \beta_kc}{p_k\tL}\right\rceil.
		\end{aligned}
	\end{align}
	For all $k\ge \Delta+1$:
	\begin{align}
		\begin{aligned}
			& \tau_k = \tau,\ \lambda_k = \lambda:=\tfrac{\tau}{1+\tau},\ \beta_k = \Delta\lambda^{-(k-\Delta)},\ p_k = \tfrac{2\tL}{1+\tau}, T_k = \left\lceil \tfrac{2(1+\tau)R\|\cA\|}{\tL\lambda^{\tfrac{k-\Delta}{2}}}\right\rceil, c_k = \left\lceil\tfrac{(1+\tau)^{2}\Delta c}{\tL^2\lambda^{\frac{k+N-2\Delta}{2}}}\right\rceil.
		\end{aligned}
	\end{align}
	And for all $k$ and $t$,
	\begin{align}
            \nonumber
		\eta_k^t &= (p_k+\mu)(t-1) + p_kT_k, \ q_k^t = \tfrac{\tL T_k}{4 \beta_kR^2}, \\
            \alpha_k^t &= \begin{cases}
			\tfrac{\beta_{k-1} T_k}{\beta_k T_{k-1}}&k\ge 2\ \text{ and } t=1
			\\
			1 & \text{otherwise}.
		\end{cases}
	\end{align}
\begin{itemize}[leftmargin=0.1in]
    \item If problem \eqref{eqn: problem} is smooth and convex, then the algorithm \ref{alg:PDS} with the LO solver \ref{alg: pds-cgs} returns an $\varepsilon$ approximation solution with a sampling complexity of $O\br{\sqrt{\tilde{L}/\varepsilon} + \sigma^2/\varepsilon^2}$, communication complexity of $O\br{\norm{\cA}/\varepsilon}$ and the linear oracle complexity is $O \left( 1/\varepsilon^2 \right)$.
    \item If problem \eqref{eqn: problem} is smooth and $\mu$-strongly convex, then the algorithm \ref{alg:PDS} with the LO solver \ref{alg: pds-cgs} returns an $\varepsilon$ approximation solution with a sampling complexity $O\br{\sqrt{\tilde{L}/\mu} \log \br{\tilde{L}/\varepsilon} + \sigma^2/\varepsilon}$, a communication complexity of $O\br{\norm{\cA}/\sqrt{\varepsilon}}$ and the linear oracle complexity is $O \left( 1/\varepsilon^2 \right)$.
\end{itemize}
\end{theorem}

Note that the gradient complexities in Theorem \ref{thm:bigtheorem} are independent of the spectral gap of the communication graph $\rho$, thus they are independent of the graph topology. Due to space limits, we only present a proof sketch here, and the complete proof is deferred to the appendix.

\begin{proof}[Proof Sketch] We will bound the primal gap $f(x_k)-f^*$ and the consensus gap $\norm{\cA x_k}$ using the following gap function with $w := (x,y,z), \hat{w}_k := (\hat{x}_k, y_k, \hat{z}_k)$:
\begin{align}
    Q(\hat{w}_k,w) := &\left[ \mu \nu(\hat{x}_k) + \langle \hat{x}_k, y + \cA^Tz \rangle - \tilde{f}^*(y) \right] - \left[ \mu \nu(x) + \langle x, y + \cA^Tz \rangle - \tilde{f}^*(y) \right].
\end{align}
The gap function $Q$ is the duality gap from the dual saddle point problem \eqref{eqn: saddle point dual} and will be used to obtain convergence guarantees on our primal-dual updates. Our proof will follow $3$ main steps:

\textbf{Step 1}: First, we will establish bounds on the gap function with the following Lemma:
\begin{lemma}
\label{pro: gap-function-bound}
Let $\varepsilon_i^t$ is the obtained error (the primal error) after running Algorithm \ref{alg: pds-cgs} when solving for $x_i^t$ and $\hat{w}_k:=(\hat{x}_k,y_k,\hat{z}_k)$, we have:
\begin{align}
    \label{eq:Qest_raw}
    &\sum_{k=1}^{N}\beta_kQ(\hat{w}_k, w) + A + B \leq C + D + \sum_{k=1}^N \frac{\beta_k}{T_k}\left(\sum_{t=1}^{T_k}\varepsilon_k^t\right) \notag\\
    &\textrm{ where } \\
    &A: = \sum_{k=1}^{N}\beta_k\left[- \langle \hx_k, y\rangle + \langle x, y_k\rangle + \langle v_k, \hx_k - x\rangle - \tf^*(v_k) + \tf^*(y) 
        	 + \frac{p_k}{T_k}\sum_{t=1}^{T_k}V(x_{k-1}, x_k^t)\right],\notag\\
    &B:= \sum_{k=1}^{N}\tfrac{\beta_k}{T_k}\sum_{t=1}^{T_k}\left[\langle \cA^\top z_k^t - \cA^\top z, x_k^t - \tu_k^t\rangle + q_k^tU(z_k^{t-1}, z_k^t)  + \eta_k^tV(x_k^{t-1}, x_k^t)\right],\notag\\
    & C:= \sum_{k=1}^{N}\tfrac{\beta_k}{T_k}\sum_{t=1}^{T_k}\left[q_k^tV(z_k^{t-1}, z) - q_k^tV(z_k^t,z)\right],\notag\\
    &D:= \sum_{k=1}^{N}\tfrac{\beta_k}{T_k}\sum_{t=1}^{T_k}\left[\eta_k^t V(x_k^{t-1}, x) - (\mu+\eta_k^t+p_k) V(x_k^t, x) + p_kV(x_{k-1}, x)\right].\notag
	\end{align}
\end{lemma}
Roughly speaking, $A$ concerns with controlling the error of the stochastic gradient oracle via the batch size, $B$ deals with the error from communication. Quantities $C, D$ control the error from the proximal functions $V$ with respect to the variables $z$ and $x$ respectively. Furthermore, Lemma \ref{pro: gap-function-bound} also provides a quantifiable relationship between the LO error and the gap functions, which allows us to control the solution accuracy using these LO errors. Additionally, the LO error can be controlled by specifying a pre-determined number of LO oracle calls at each inner iteration. This step will be a stepping stone to establishing bounds on the primal gap and consensus gap below.
    
\textbf{Step 2}: Next, we will convert the bound in step $1$ to obtain bounds on the primal gap $f(x_N)-f^*$ and the consensus gap $\norm{\cA x_N}$ to arrive at the following Proposition.


\begin{prop}
    \label{pro: LO-error-proposition-stochastic}
    Let $\varepsilon_k^i$ be the error from running the LO oracle at the $i$-th inner iteration and $k$-th outer iteration. Suppose that the parameters have been chosen as in Theorem \ref{thm:bigtheorem}. Under Assumptions \ref{assumption: agent-smoothness} \ref{assumption: stochastic-gradient-assumptions}, we have:


	\begin{align}
    	& \begin{aligned}
    	& \E[f(x_N) - f(x^*)]
    		\leq  \br{\sum_{k=1}^{N}\beta_k}^{-1} \left( \sum_{k=1}^N \frac{\beta_k}{T_k}\left(\sum_{t=1}^{T_k} \varepsilon_k^t\right) + \beta_1\left(\frac{\eta_1^1}{T_1} + p_1\right)V(x_0, x^*) + \sum_{k=1}^N \frac{\beta_k \sigma}{p_k c_k}\right),
    	\end{aligned}
\\
    		& \begin{aligned}
               & \E[\|\cA x_N\|] \leq \br{\sum_{k=1}^{N}\beta_k}^{-1} \left( \sum_{k=1}^N \br{\frac{\beta_k}{T_k}\left(\sum_{t=1}^{T_k} \varepsilon_k^t\right) + \frac{\beta_k \sigma}{p_k c_k}} + \beta_1 \left[\frac{q_1^1}{2T_1}(\|z^*\|_2+1)^2 + \left(\tfrac{\eta_1^1}{T_1} + p_1\right)V(x_0, x^*)\right] \right)
    	\end{aligned}
	\end{align}	
\end{prop}
The full description of this Proposition is reserved to the Appendix. In essence, Proposition \ref{pro: LO-error-proposition-stochastic} provides us a way to control the LO error by bounding the mean primal gap and the mean consensus gap with the LO error, which will eventually pave the way to analyze the complexity of the algorithm. To prove this Proposition, we utilize the convexity of the gap function $Q$. Let $\hat{w}_k = (\hat{x}_k, y_k, \hat{z}_k)$, the convexity property of $Q$ implies:
\begin{align}
    \frac{\sum_{t=1}^{T_k} \varepsilon_k^t}{T_k} \geq \frac{\sum_{t=1}^{T_k} Q(w_k^t,w)}{T_k} \geq Q(\hat{w}_k,w) \, \text{ and } \\
    \sum_{k=1}^N\beta_kQ(\hat{w}_k,w) \geq (\sum_{k=1}^N\beta_k)Q(\overline{w}_N,w).
\end{align}
From the Lemma \ref{pro: gap-function-bound} and from the observation that $Q(\overline{w}_N, w) \geq \E[f(\overline{x}_N)-f^*]$ and $Q(\overline{w}_N, \overline{w}_N^*) \geq \E[\norm{\cA x_N}]$, we further bound the quantities $A,B,C,D$ in Lemma \ref{pro: gap-function-bound} to arrive at the primal gap and consensus gap bounds. Note that the quantities $C,D$ can be bound using a telescoping sum and thus are upper bounded by some initial terms. On the other hand, the quantity $-A$ can be upper-bounded by the weighted sum of the variance of the stochastic gradient oracle and the quantity $B$ is non-negative, which means that $-A-B$ is upper-bounded by the weighted sum of the variance of the stochastic gradient oracle.

\textbf{Step 3}: Now, toward the last step, we will choose suitable parameters to obtain competitive finite-time bounds from Proposition \ref{pro: LO-error-proposition-stochastic}. To establish the linear oracle complexity bounds, we will choose the number of LO calls at each vertex and at specific outer iterations such that the total error incurred by the linear oracles is at most $\varepsilon/2$. At the same time, we choose a suitable number of outer iteration $N$ such that the primal error and the consensus error without the linear oracle error is also at most $\varepsilon/2$. Since $\varepsilon$ and $\varepsilon/2$ only differ by a constant, the outer and inner iteration complexity is still the same with or without the linear oracle error. 
\end{proof}

\textbf{Remark}: In the case of using a full gradient to update the algorithm (which corresponds to when $\sigma = 0$), we have the gradient complexity is $O\br{1/\sqrt{\varepsilon}}$, the communication complexity is $O\br{1/\varepsilon}$ and the LO complexity is $O\br{1/\varepsilon^2}$ for when $f$ is smooth and convex. If $f$ is smooth and $\mu$ strongly-convex, we have the gradient complexity is $O\br{\log 1/\varepsilon}$, the communication complexity is $O\br{1/\sqrt{\varepsilon}}$ and the LO complexity is $O\br{1/\varepsilon^2}$. In such cases, we still have the gradient complexity better than that of \cite{wai-fw-2017}, in addition to not having to rely on the solution lying strictly inside the constraint set. Such an assumption is not practical since many constraints like the $\ell_1$ norm constraint will have the solution lying on the boundary of the constraint set.

\section{Numerical experiments}
\label{sec:experiments}
\subsection{Logistics regression experiments}
\begin{figure}[ht!]
\vspace{-0.2in}
  \centering
  \subfigure[Primal gap of each method for $\mu = 0$]{\includegraphics[scale=0.49]{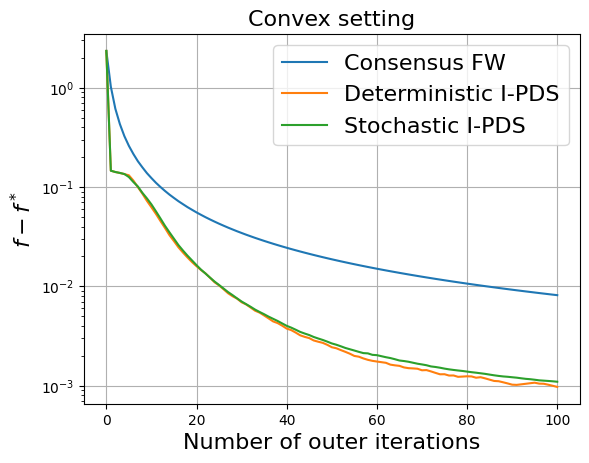}}
  \subfigure[Primal gap of each method for $\mu = 0.5$]{\includegraphics[scale=0.49]{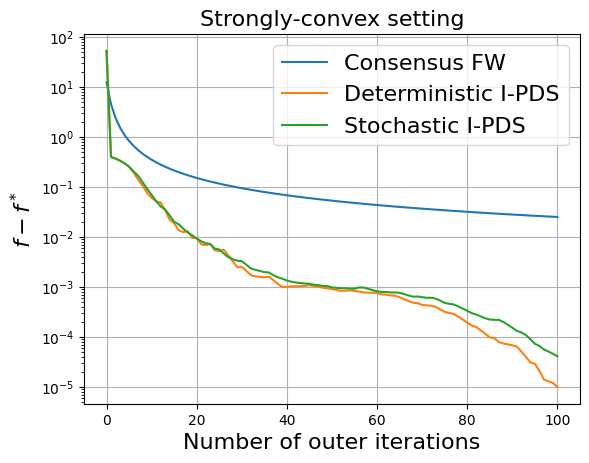}} \\
  \subfigure[LO complexity of each method for $\mu = 0$]{\includegraphics[scale=0.49]{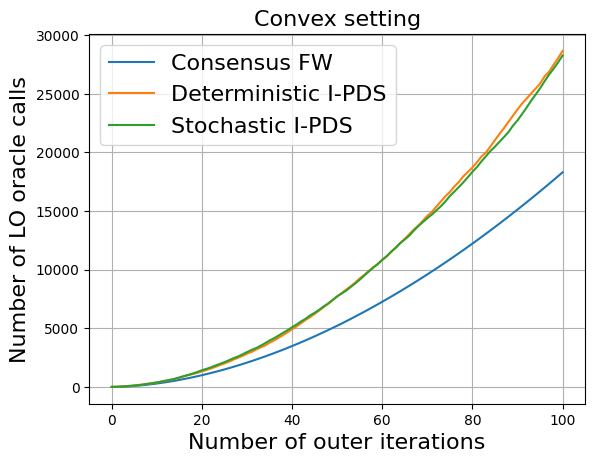}}
  \subfigure[LO complexity of each method for $\mu = 0.5$]{\includegraphics[scale=0.49]{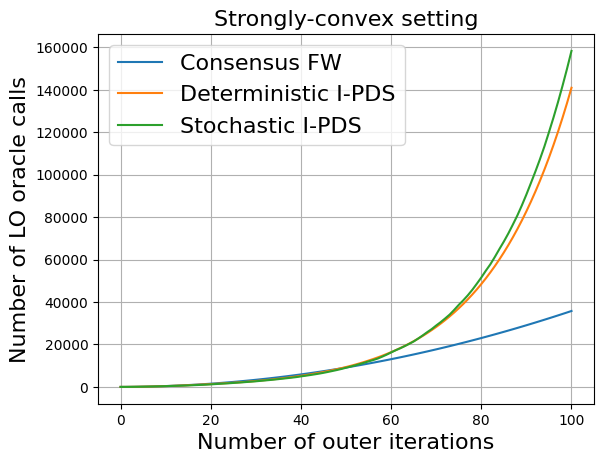}}
\end{figure}


In this section, we demonstrate the advantages of our proposed I-PDS method through some preliminary numerical experiments and compare it with the projected gradient method in addition to the consensus-based methods proposed by \cite{wai-fw-2017}. We also use the oracle scheme in \cite{nemirovski-cg-2013}, \cite{complexity-LO-oracles} for our setting. We consider a decentralized convex smooth optimization problem of the unregularized logistic regression model. The dataset of our choice will be the ijcnn1 datasets obtained from LIBSVM, which is not linearly separable, and we choose $20000$ samples from this dataset. For the linearly constrained problem \ref{eqn: primal-dual problem}, we choose $\cA = \cL \otimes I_d$ where $\cL$ is the Laplacian matrix of the graph $G$, whose entries are $|N_i|-1$ for any $i = j$ where $|N_i|$ is the degree of the vertex $i$ and $-1$ for any $i \neq j, (i,j) \in \cE$. We also generate the communication network using the Erdos-Renyi algorithm. Our generated graph will have $m = 10$ nodes and we will split the $20000$ samples of our dataset over these nodes evenly. Initially, all nodes will have the same initial point $x_0 = y_0 = z_0 = 0$. We then compare the performance of the I-PDS algorithm with the consensus-based decentralized Frank-Wolfe algorithm DeFW by \cite{wai-fw-2017}. For the stochastic I-PDS (when $\sigma > 0$, denoted as SPDS in the plots), since there is no comparable counterpart, we will simply report its performance. At each outer iteration, we sample $m$ data points from the dataset and distribute them over $m$ nodes with $1$ each and then update the algorithm with the stochastic gradient estimate. We report our results in the following figures. We note that the wall-clock time of the LO oracle is highly dependent on the implementation and our main interest in this work is to reduce the number of data oracle access. 



In addition to the primal gap losses, we also report the number of gradient samples after running $100$ outer iterations on a dataset with $20000$ data points of each method in Table \ref{tab:data-oracle-calls}. Observe that the Stochastic PDS method has a significantly smaller number of gradient samples for the same number of outer iterations.

\begin{table}[H]
\begin{center}
\small
\begin{tabular}{|c|c|c|}
    \hline
    Method & Convex & Strongly-convex \\
    \hline
DeFW & $2 \times 10^7$ & $2 \times 10^7$      \\
\hline
Deterministic I-PDS $(\sigma = 0)$ & $2 \times 10^7$ & $2 \times 10^7$ \\ \hline
Stochastic I-PDS $(\sigma > 0)$ & $21200$ & $78740$ \\ \hline
\end{tabular}
\caption{Comparison of the DeFW and PDS algorithms (deterministic and stochastic) in terms of the number of gradient samples after $100$ outer iterations. \label{tab:data-oracle-calls}}
\end{center}
\end{table}

\subsection{The effects of graph topologies}
\label{sec:additional-experiments}
To empirically validate the graph topology independence property of our algorithm and the discussions in Section \ref{sec:spectral-gap}, we will compute the primal gap per the number of gradient evaluations of each algorithm for different graph topologies. We follow the same setting in Section \ref{sec:experiments} and measure the number of gradient evaluations required to reach the target loss of $70$ in different graph topologies with $100$ vertices each. While DeFW requires more than $1000000$ gradient evaluations for all graph topologies, we note that for any fixed number of gradient evaluations, DeFW achieves the smallest primal gap in the complete graph and the largest in the path graph. This confirms the graph topology dependence of DeFW where DeFW converges slower for graphs with smaller spectral gaps. In contrast, I-PDS requires roughly the same amount of gradient evaluations to reach the target loss.

\begin{table}[H]
\begin{center}
\small
\begin{tabular}{|c|c|c|c|}
    \hline
    Method & Graph & Loss & \begin{tabular}[c]{@{}c@{}}Gradient \\ evaluations\end{tabular} \\
    \hline
DeFW & Erdos-Renyi $(p = 0.1)$ & $70$                      & $\geq 10^6$        \\
\hline
I-PDS  &  Erdos-Renyi $(p = 0.1)$  & $70$   & $866$   \\ \hline
DeFW & Path Graph     & $70$ & $\geq 10^6$  \\ \hline
I-PDS  & Path Graph    & $70$  & $865$    \\ \hline
DeFW  & Barbell Graph &$70$     & $\geq 10^6$                                  \\ \hline
I-PDS & Barbell Graph & $70$ & 865               \\ \hline 
DeFW  & Complete Graph &$70$     & $\geq 10^6$                                  \\ \hline
I-PDS & Complete Graph & $70$ & $865$               \\ \hline 
\end{tabular}
\caption{Comparison of the DeFW and PDS algorithms in terms of reaching the same target loss for different graph topologies with $100$ vertices each. \label{tab:different-topologies}}
\end{center}
\end{table}

\section{Discussion}
\label{sec: discussion}
In this paper, we have studied a graph topology invariant decentralized projection-free algorithm for constrained optimization problems. With the I-PDS framework, we require fewer data oracle access and match the communication complexity with the best-known results while allowing an inexact gradient and the solution to be on the boundary. Furthermore, our gradient sampling complexity is graph topology invariant and will not penalize networks with small spectral gaps. We are aware that the LO complexity of $O \br{1/\varepsilon^2}$ is worse than the existing results and we suspect that improving the LO complexity given the current assumptions will be difficult since the LO complexity of CGS matches the lower bound of \cite{guzman2018lower}. In future works, we might be able to get improved rates on the LO complexity in different settings where faster rates are possible and apply our method to a variety of applications such as Computational Optimal Transport \cite{quangnguyen2022unbalanced,ducnguyen2023partial-ot}, Semidefinite Programming \cite{arnesh2024lowrank-sdp} or Stochastic Approximation \cite{hoang_nonlinear_sa}.

\bibliographystyle{ieeetr}
\bibliography{refs}

\clearpage
\onecolumn
\appendix

\section{Discussions on the spectral gap of the communication network}
\label{sec:spectral-gap}
Recall that from Table \ref{complexity-table}, our gradient complexity is better when $M\rho^{-1} >> \sigma^2/\varepsilon$ in the convex case and $M\rho^{-1} >> \sigma^2/\sqrt{\varepsilon}$ where $M$ is the number of training data points. Thus, it is critical to discuss the size of the spectral gap in different graph topologies, given that the gradient sampling complexities of our method are graph topology invariant. We present several graph topologies with the corresponding spectral gap below (with $m$ being the number of vertices):

\begin{table}[htb!]
\centering
\label{spectral-gap-table}
\begin{tabular}{|c | c|} 
 \hline
 Graph topology & Spectral gap \\
 [0.5ex] 
 \hline
 \multirow{2}{*}{Chain graph} & \multirow{2}{*}{$O\br{1/m^2}$} \\
 & \\
 \hline
 \multirow{2}{*}{Star graph} & \multirow{2}{*}{$O\br{1/m^2}$} \\
 & \\
 \hline 
 \multirow{2}{*}{Geometric random graph} & \multirow{2}{*}{$O\br{1/(m \log m)}$} \\
 & \\
 \hline
 \multirow{2}{*}{Barbell graph} & \multirow{2}{*}{$O\br{1/m^3}$} \\
 & \\
 \hline
 \multirow{2}{*}{Cubic graph} & \multirow{2}{*}{$O\br{1/m^2}$} \\
 & \\
 \hline
\end{tabular}
\label{table}
\caption{ Common graph topologies with a small spectral gap \cite{nedic2018network, aldous-random-walk-book}}
\end{table}
Under these topologies (which could be very common in many applications), the gradient sampling complexity of the consensus-based algorithm could be very large as it scales with the number of workers $m$. This could severely hinder the training of large-scale systems. In contrast to consensus-based algorithms in which the unfavorable dependence of spectral gap is inevitable, our primal-dual sliding approach is graph topology invariant.

\section{Proof of key results}
In the following section, we will present the missing proofs of the results presented in Section \ref{sec:main-reults}.

\subsection{Linear oracle error bounds}
In order to establish bounds on the linear oracle error, we need to first bound the gap function:
\begin{align}
    Q(\hat{w}_k,w) := \left[ \mu \nu(\hat{x}_k) + \langle \hat{x}_k, y + \cA^Tz \rangle - \tilde{f}^*(y) \right] - \left[ \mu \nu(x) + \langle x, y + \cA^Tz \rangle - \tilde{f}^*(y) \right].
\end{align}
We have the following Proposition:
\begin{lemma}
	Suppose that $\hat{x}_k=\sum_{t=1}^{T_k}x_k^t/T_k$ and $\hz_k=\sum_{t=1}^{T_k}z_k^t/T_k$, where the iterates $\{x_k^t\}_{t=1}^{T_k}$ and $\{z_k^t\}_{t=1}^{T_k}$ are defined by 
	\begin{align}
    	\label{eq:zkt_proof}
    	z_k^t =& \argmin_{z \in \R^{md}}h(z) + \langle -\cA\tu_k^t, z\rangle + q_k^tU(z_k^{t-1}, z),
    	\\
    	\label{eq:xkt_proof}
    	x_k^t =& \argmin_{x\in \cX} \mu\nu(x) + \langle v_k + \cA^\top z_k^t, x\rangle + \eta_k^tV(x_{k}^{t-1}, x) + p_k V(x_{k-1}, x)
	\end{align}
	respectively. In addition, let $\varepsilon_i^t$ is the obtained error (the primal error) after running Algorithm \ref{alg: pds-cgs} when solving for $x_i^t$. Letting $\hat{w}_k:=(\hat{x}_k,y_k,\hat{z}_k)$ we have
	\begin{align}
    	\label{eq:Qest_raw}
    	&\sum_{k=1}^{N}\beta_kQ(\hat{w}_k, w) + A + B \leq C + D + \sum_{k=1}^N \frac{\beta_k}{T_k}\left(\sum_{t=1}^{T_k} \varepsilon_k^t\right) \forall w:=(x,y,z) \in \cX \times \R^{md} \times \R^{md}
	\end{align}
	where
	\begingroup
	\allowdisplaybreaks
	\begin{align}
    	\label{eq:A}
    	& 
                \begin{aligned}
        	A: = & \sum_{k=1}^{N}\beta_k\left[- \langle \hx_k, y\rangle + \langle x, y_k\rangle + \langle v_k, \hx_k - x\rangle - \tf^*(v_k) + \tf^*(y) 
        	 + \frac{p_k}{T_k}\sum_{t=1}^{T_k}V(x_{k-1}, x_k^t)\right],
        	\end{aligned}
    	\\
    	\label{eq:B}
    	& 
                \begin{aligned}
        	B:= & \sum_{k=1}^{N}\tfrac{\beta_k}{T_k}\sum_{t=1}^{T_k}\left[\langle \cA^\top z_k^t - \cA^\top z, x_k^t - \tu_k^t\rangle + q_k^tU(z_k^{t-1}, z_k^t)  + \eta_k^tV(x_k^{t-1}, x_k^t)\right],
        	\end{aligned}
    	\\
    	\label{eq:C}
    	& C:= \sum_{k=1}^{N}\tfrac{\beta_k}{T_k}\sum_{t=1}^{T_k}\left[q_k^tU(z_k^{t-1}, z) - q_k^tU(z_k^t,z)\right],
    	\\
    	\label{eq:D}
    	& D:= \sum_{k=1}^{N}\tfrac{\beta_k}{T_k}\sum_{t=1}^{T_k}\left[\eta_k^t V(x_k^{t-1}, x) - (\mu+\eta_k^t+p_k)V(x_k^t, x) + p_kV(x_{k-1}, x)\right].
	\end{align}
	\endgroup
\end{lemma}
\begin{proof}
From $(5.2)$ in \cite{lan2021graph}, the optimality condition gives us:
\begin{align}
    \langle (\mu + p_k + \eta_k^t)\nu'(x_k^t) + v_k + \mathcal{A}^Tz_k^t - \eta_k^t \nu'(x_k^{t-1}) - p_k\nu'(x_{k-1}), x_k^t - x \rangle \leq 0.
\end{align}
Thus, assume that the Frank-Wolfe algorithm was able to obtain an $\varepsilon_k^t$-approximation solution to the problem \eqref{linear-problem} after $s_k^t$ steps to solve for $x_k^t$, at each step it will incur an error of $\varepsilon_k^t$. This gives:
\begin{align}
    \label{eqn: LO-approximation-implication}
    \langle (\mu + p_k + \eta_k^t)\nu'(x_k^t) &+ v_k + \mathcal{A}^Tz_k^t - \eta_k^t \nu'(x_k^{t-1}) - p_k\nu'(x_{k-1}), x_k^t - x \rangle \leq \varepsilon_k^t.
\end{align}
From here, we will adapt the analysis in Proposition 5.1 in \cite{lan2021graph} as follows: from \ref{eqn: LO-approximation-implication} and from the convexity of $h$ and $\nu$, and the definitions of $U$ and $V$, we obtain the following two relations:
\begin{align}
	\label{eq:oc_zkt}
        \nonumber
	& h(z_k^t) - h(z) + \langle -\cA\tu_k^t, z_k^t - z\rangle + q_k^tU(z_k^{t-1}, z_k^t) \\
    &+ q_k^tU(z_k^t,z) \le q_k^tU(z_k^{t-1}, z), \forall z \in \R^{md},
	\\
	& 
	\begin{aligned}
	&\langle v_k+\cA^\top z_k^t, x_k^t - x\rangle + \mu\nu(x_k^t) - \mu\nu(x) + \eta_k^tV(x_k^{t-1}, x_k^t) \\
        &+ (\mu + \eta_k^t+p_k)V(x_k^t, x)  + p_kV(x_{k-1}, x_k^t)
	\\
	&\le \eta_k^t V(x_k^{t-1}, x) + p_kV(x_{k-1}, x) + \varepsilon_k^t, \forall x \in \cX.
	\end{aligned}
	\label{eq:oc_xkt}
\end{align}
Notice that the second inequality now has an additional $\varepsilon_k^t$ on the RHS rather than $0$ as in the original analysis. Summing up the two relations above and notice that with the identity:
\begin{align}
        \nonumber
	\langle -\cA\tu_k^t, z_k^t - z\rangle &+ \langle v_k+\cA^\top z_k^t, x_k^t - x\rangle 
	\\
        \nonumber
	= & \langle \cA^\top z_k^t - \cA^\top z, x_k^t - \tu_k^t\rangle + \langle \cA^\top z, x_k^t\rangle \\
        - & \langle v_k+\cA^\top z_k^t, x\rangle + \langle v_k, x_k^t\rangle,
\end{align}
we have:
\begin{align}
        \nonumber
	& \langle \cA^\top z_k^t - \cA^\top z, x_k^t - \tu_k^t\rangle + \langle \cA^\top z, x_k^t\rangle - \langle v_k+\cA^\top z_k^t, x\rangle \\
        \nonumber
        & + \langle v_k, x_k^t\rangle + h(z_k^t) - h(z) +  q_k^tU(z_k^{t-1}, z_k^t) \\
        \nonumber
        & + q_k^tU(z_k^t,z) + \mu\nu(x_k^t) - \mu\nu(x) + \eta_k^tV(x_k^{t-1}, x_k^t) \\
        \nonumber
        & + (\mu + \eta_k^t+p_k)V(x_k^t, x)  + p_kV(x_{k-1}, x_k^t)
	\\
	& \le q_k^tU(z_k^{t-1}, z) + \eta_k^t V(x_k^{t-1}, x) + p_kV(x_{k-1}, x) + \varepsilon_k^t.
\end{align}
Summing from $t=1,\ldots,T_k$ and noting the definitions of $\hx_k$ and $\hz_k$ and the convexity of functions $h$ and $\nu$ we have:
\begin{align}
        \nonumber
	& \ T_k \left[\langle \cA^\top z, \hx_k\rangle - \langle v_k+\cA^\top \hz_k, x\rangle + \langle v_k, \hx_k\rangle + h(\hz_k) - h(z) + \mu\nu(\hx_k) - \mu\nu(x)\right] \\
        \nonumber
        &+ \sum_{t=1}^{T_k}p_k \left[V(x_{k-1}, x_k^t) + \langle \cA^\top z_k^t - \cA^\top z, x_k^t - \tu_k^t\rangle + q_k^tU(z_k^{t-1}, z_k^t)  + \eta_k^tV(x_k^{t-1}, x_k^t)\right]  \\
 &\le  \sum_{t=1}^{T_k}\left[q_k^tU(z_k^{t-1}, z) - q_k^tU(z_k^t,z) + \eta_k^t V(x_k^{t-1}, x) - (\mu + \eta_k^t+p_k)V(x_k^t, x) + p_kV(x_{k-1}, x) + \varepsilon_k^t\right].
\end{align} 
Multiplying the last inequality by $\beta_k/T_k$ and noting the definition of the gap function $Q(\hat{w}_k,w)$ where $\hat{x}_k =  \frac{\sum_{t=1}^{T_k} x_k^t}{T_k}, \hat{z}_k = \frac{\sum_{t=1}^{T_k} z_k^t}{T_k}$, we obtain the bound:
\begin{align}
    &\sum_{k=1}^N \frac{\beta_k}{T_k} \br{\sum_{t=1}^{T_k} \varepsilon_k^t} + C + D \geq A + B + \sum_{k=1}^N \beta_k Q(\hat{w}_k, w), \forall w \in \cX \times \R^{md} \times \R^{md}.
\end{align}
Hence proved.
\end{proof}
From the established bounds on the gap function, we can subsequently bound the quantities $A,B,C,D$ to obtain the bounds on the primal gap and the consensus gap in relation to the preset parameters. 

\subsubsection{Proof of Proposition \ref{pro: LO-error-proposition-stochastic}}
\begin{prop}
	Suppose that Assumption \ref{assumption: agent-smoothness} is satisfied and assume that the parameters of Algorithm \ref{alg:PDS} satisfy the following conditions:

		$\bullet$ For any $k\ge 2$, 
		\begin{align}
    			\label{eq:cond_Qest_k}
    			\begin{aligned}
    			& \beta_k\tau_k\le \beta_{k-1} (\tau_{k-1}+1),\ \beta_{k-1} = \beta_k\lambda_k,\tL\lambda_k\le p_{k-1}\tau_k, \beta_k T_{k-1}\alpha_k^1= \beta_{k-1} T_k, \\
                    &\alpha_k^1 \|\cA\|^2\le \eta_{k-1}^{T_{k-1}}q_k^{1}, \beta_k T_{k-1} q_k^1\le \beta_{k-1} T_k q_{k-1}^{T_{k-1}},
    			\\
    			& \beta_k T_{k-1} (\eta_k^1 + p_k T_k)\le \beta_{k-1} T_k(\mu+\eta_{k-1}^{T_{k-1}}+p_{k-1});
    		    \end{aligned}
		\end{align}
		$\bullet$  For any $t\ge 2$ and $k\ge 1$, 
		\begin{align}
			\label{eq:cond_Qest_t}
			\begin{aligned}
				& \alpha_k^t=1,\ \|\cA\|^2\le\eta_k^{t-1} q_k^t, \ q_k^t \le q_k^{t-1}, \eta_k^t\le \mu+\eta_k^{t-1}+p_k;
			\end{aligned}
		\end{align}
		$\bullet$  In the first and last outer iterations, 
		\begin{align}
			\label{eq:cond_Qest_N}
			\begin{aligned}
				& \tau_1 = 0,\ p_N(\tau_N+1)\ge \tL, \text{ and }
				\eta_N^{T_N}q_{N}^{T_N}\ge \|\cA\|^2.
			\end{aligned}
		\end{align}

	Also, let $\varepsilon_i^t$ is the obtained error (the primal error) after letting Algorithm \ref{alg: pds-cgs} runs for $w_i^t$ iterations. Then we have
	\begin{align}
    	& \begin{aligned}
    	& \E[f(x_N) - f(x^*)]
    		\leq  \br{\sum_{k=1}^{N}\beta_k}^{-1} \left( \sum_{k=1}^N \frac{\beta_k}{T_k}\left(\sum_{t=1}^{T_k} \varepsilon_k^t\right) + \beta_1\left(\frac{\eta_1^1}{T_1} + p_1\right)V(x_0, x^*) + \sum_{k=1}^N \frac{\beta_k \sigma}{p_k c_k}\right),
    	\end{aligned}
\\
    		& \begin{aligned}
    		& \E[\|\cA x_N\|_2] \leq 
    			 \frac{\sum_{k=1}^N \br{\frac{\beta_k}{T_k}\left(\sum_{t=1}^{T_k} \varepsilon_k^t\right) + \frac{\beta_k \sigma}{p_k c_k}} + \beta_1
    			 \left[
    			 \frac{q_1^1}{2T_1}(\|z^*\|_2+1)^2 + \left(\tfrac{\eta_1^1}{T_1} + p_1\right)V(x_0, x^*)\right]}{\sum_{k=1}^N \beta_k}.
    	\end{aligned}
	\end{align}	
\end{prop}

\begin{proof}
Similar to the proof of Proposition \ref{pro: gap-function-bound} and Proposition 3.1 in \cite{lan2021graph}, we will study the gap function $Q(\cdot,\cdot)$. We have:
\begin{align}
    	\label{eq:Qest_raw}
    	\sum_{k=1}^{N}\beta_kQ(\hat{w}_k, w) + A + B \leq C + D + \sum_{k=1}^N \frac{\beta_k}{T_k}\left(\sum_{t=1}^{T_k} \varepsilon_k^t\right)
\end{align}
$\forall w:=(x,y,z)\in\cX \times \R^{md} \times \R^{md}$. Recall that we have the relations: $\hat{x}_k =  \frac{\sum_{t=1}^{T_k} x_k^t}{T_k}, \hat{z}_k = \frac{\sum_{t=1}^{T_k} z_k^t}{T_k}$. By the convexity of $Q$, we have for all $w \in \cX \times \cY \times \cZ$:
\begin{align}
    \frac{\sum_{t=1}^{T_k} \varepsilon_k^t}{T_k} \geq \frac{\sum_{t=1}^{T_k} Q(w_k^t,w)}{T_k} \geq Q(\hat{w}_k,w).
\end{align}
and
\begin{align}
    \sum_{k=1}^N\beta_kQ(\hat{w}_k,w) \geq (\sum_{k=1}^N\beta_k)Q(\overline{w}_N,w).
\end{align}
Thus, what is left to be done is to obtain the bounds on the quantities $A, B, C, D$. From the inequality $Q(\overline{w}_N, w) \geq \E[f(\overline{x}_N)-f^*]$ and apply lemmas $5.3, 5.4, 5.5$ in \cite{lan2021graph}, we obtain the bound for the primal gap as follows:
\begin{align}
    \label{eqn: stochastic-primal-gap}
    \E[f(x_N) - f(x^*)]
    		\leq \frac{\beta_1\left(\frac{\eta_1^1}{T_1} + p_1\right)V(x_0, x^*) + \sum_{k=1}^N \br{\frac{\beta_k \sigma}{p_k c_k} + \frac{\beta_k}{T_k}\left(\sum_{t=1}^{T_k} \varepsilon_k^t\right)}}{\sum_{k=1}^{N}\beta_k}.
\end{align}
And from the inequality $Q(\overline{w}_N, \overline{w}_N^*) \geq \E[\norm{\cA x}]$, the bound for consensus gap can be obtained as follows:
\begin{align}
    \label{eqn: stochastic-consensus-gap}
    &\E[\|\cA x_N\|_2] \leq \frac{\beta_1
    			 \left[
    			 \frac{q_1^1}{2T_1}(\|z^*\|_2+1)^2 + \left(\tfrac{\eta_1^1}{T_1} + p_1\right)V(x_0, x^*)\right] + \sum_{k=1}^N \br{\frac{\beta_k \sigma}{p_k c_k} + \frac{\beta_k}{T_k}\left(\sum_{t=1}^{T_k} \varepsilon_k^t\right)}}{\sum_{k=1}^{N}\beta_k}.
\end{align}
Hence proved.
\end{proof}

\subsection{Main results}
In this Section, we will present the proof of Theorem \ref{thm:bigtheorem}, which consists of two main cases: the convex setting and the strongly-convex setting.

First, we will present the following Lemma on the guarantee of the CGS procedure
\begin{lemma}
Algorithm \ref{alg: pds-cgs} produces an $\varepsilon_k^t$-approximation solution within $\left\lfloor \frac{12(\eta_k^t+p_k)D_{\cX}^2}{\varepsilon_k^t} \right\rfloor$
\end{lemma}
The proof is this Lemma can be found in \cite{Lan2016ConditionalGS}. Thus, the LO complexity of our algorithm will be heavily dependent on how we choose our LO error tolerance. Now, return to the main Theorem, we will present the proof below.

\begin{theorem}
Denote $N$ as the pre-determined number of outer iterations, $\tau:=2\sqrt{\tL/\mu}$ and $\Delta:=\lceil 2\tau + 1\rceil$ if $\mu>0$, and $\Delta:=+\infty$ if $\mu=0$. Suppose that the Assumptions \ref{assumption: agent-smoothness}, \ref{assumption: stochastic-gradient-assumptions} hold, and that the parameters in Algorithm \ref{alg:PDS} are set to the following: 

        For all $k\ \leq \Delta$:
	\begin{align}
		\label{eq:par_beforeDelta_S}
		\begin{aligned}
            \nonumber
			&\tau_k = \tfrac{k-1}{2},\ \lambda_k = \tfrac{k-1}{k},\ \beta_k = k,\ p_k = \tfrac{4\tL}{k}, T_k = \left\lceil \tfrac{kR\|\cA\|}{\tL}\right\rceil, c_k = \left\lceil \tfrac{\min\{N,\Delta\} \beta_kc}{p_k\tL}\right\rceil.
		\end{aligned}
	\end{align}
	For all $k\ge \Delta+1$:
	\begin{align}
		\begin{aligned}
			& \tau_k = \tau,\ \lambda_k = \lambda:=\tfrac{\tau}{1+\tau},\ \beta_k = \Delta\lambda^{-(k-\Delta)},\ p_k = \tfrac{2\tL}{1+\tau},\ 
			\\
			& T_k = \left\lceil \tfrac{2(1+\tau)R\|\cA\|}{\tL\lambda^{\tfrac{k-\Delta}{2}}}\right\rceil, c_k = \left\lceil\tfrac{(1+\tau)^{2}\Delta c}{\tL^2\lambda^{\frac{k+N-2\Delta}{2}}}\right\rceil.
		\end{aligned}
	\end{align}
	And for all $k$ and $t$,
	\begin{align}
            \nonumber
		& \eta_k^t = (p_k+\mu)(t-1) + p_kT_k, \ q_k^t = \tfrac{\tL T_k}{4 \beta_kR^2}, \\
            & \alpha_k^t = \begin{cases}
			\tfrac{\beta_{k-1} T_k}{\beta_k T_{k-1}}&k\ge 2\ \text{ and } t=1
			\\
			1 & \text{otherwise}.
		\end{cases}
	\end{align}
Let $V(\cdot,\cdot) = \norm{\cdot-\cdot}_2^2$. If problem \eqref{eqn: problem} is smooth and convex, then the algorithm \ref{alg:PDS} with the LO solver \ref{alg: pds-cgs} returns an $\varepsilon$ approximation solution with a sampling complexity of $O\br{\frac{\tilde{L}}{\sqrt{\varepsilon}} + \frac{\sigma^2}{\varepsilon^2}}$ and communication complexity of $O\br{\frac{1}{\varepsilon}}$. Otherwise, if problem \ref{eqn: problem} is smooth and $\mu$-strongly convex, then the algorithm \ref{alg:PDS} with the LO solver \ref{alg: pds-cgs} returns an $\varepsilon$ approximation solution with a sampling complexity $O\br{\log \frac{1}{\varepsilon} + \frac{\sigma^2}{\varepsilon}}$ and a communication complexity of $O\br{\frac{1}{\sqrt{\varepsilon}}}$. In addition, we can obtain an $\varepsilon$-solution with the linear oracle complexity for both settings is $O \left( \frac{1}{\varepsilon^2} \right)$.
\end{theorem}

\begin{proof}
\textbf{Convex setting}: In this setting, we have $\mu = 0$ hence we will always have $k \leq \Delta$ as we have $\Delta = +\infty$. This implies that the value of our parameters will be $\eta_k^t = (p_k + \mu)(t-1) + p_kT_k, p_k = \frac{2L}{k}, T_k = \left \lceil \frac{kR\norm{\mathcal{A}}}{L} \right \rceil, \beta_k = k$. Denote $\varepsilon_k^t, s^t_k$ as the obtained error and the number of LO iterations corresponds to the value at node $t$ and in the $k$-th inner iteration, we have from \cite{Lan2016ConditionalGS} that with the $2(\eta_k^t+p_k)$-smooth LO objective and diameter $D_{\cX}$ of the constraint set $\cX$, we have:
\begin{align}
    \varepsilon_k^t \leq \frac{12(\eta_k^t+p_k)D_{\cX}^2}{s_k^t+1}
\end{align}
Choose $s_k^t = \left \lfloor \frac{24(\eta_k^t+p_k)D_{\cX}^2}{\varepsilon} \right \rfloor$, we will show that this choice of the number of LO iterations will yield an $\varepsilon$-approximation solution. Indeed, we have:
\begin{align}
     \frac{\sum_{k=1}^N \frac{\beta_k}{T_k}(\sum_{t=1}^{T_k} \epsilon_k^t)}{\sum_{k=1}^N \beta_k} \leq \frac{\sum_{k=1}^N \frac{\beta_k}{T_k}\br{\sum_{t=1}^{T_k} \frac{12(\eta_k^t+p_k)D_{\cX}^2}{s_k^t+1}}}{\sum_{k=1}^N \beta_k} \leq \frac{\sum_{k=1}^N \frac{\beta_k}{T_k}\br{\sum_{t=1}^{T_k} \frac{12(\eta_k^t+p_k)D_{\cX}^2}{\frac{24(\eta_k^t+p_k)D_{\cX}^2}{\varepsilon}}}}{\sum_{k=1}^N \beta_k} = \frac{\varepsilon}{2}
\end{align}
where $\varepsilon_k^t$ is the obtained LO error using the Frank-Wolfe method. Note that from the choice of our parameters, we have that:
\begin{align}
		\frac{\beta_1\left(\frac{\eta_1^1}{T_1} + p_1\right)V(x_0, x^*) + \sum_{k=1}^N \frac{\beta_k \sigma}{p_k c_k}}{\sum_{k=1}^{N}\beta_k} &\le \frac{2}{N^2}\left[ 4\tL \|x_0 - x^*\|_2^2 + \tfrac{\tL\sigma^2}{c}\right],
		\\
            \nonumber
		\frac{\beta_1
    			 \left[
    			 \frac{q_1^1}{2T_1}(\|z^*\|_2+1)^2 + \left(\tfrac{\eta_1^1}{T_1} + p_1\right)V(x_0, x^*)\right] + \sum_{k=1}^N \frac{\beta_k \sigma}{p_k c_k}}{\sum_{k=1}^{N}\beta_k} &\le \frac{2}{N^2} \left[\tfrac{\tL}{8 R^2}(\|z^*\|_2+1)^2 \right. \\
                &\quad \left.+ 4\tL \|x_0 - x^*\|_2^2 + \tfrac{\tL\sigma^2}{c}\right].
	\end{align}
Choose $N = \left\lceil \sqrt{{40 \tilde{L} \|x_0 - x^*\|_2^2}/{\varepsilon}}\right\rceil$, we have that $N = O\br{\sqrt{\frac{\tilde{L}}{\varepsilon}}}$ such that:
\begin{align}
    \frac{\beta_1\left(\frac{\eta_1^1}{T_1} + p_1\right)V(x_0, x^*) + \sum_{k=1}^N \frac{\beta_k \sigma}{p_k c_k}}{\sum_{k=1}^{N}\beta_k} &\leq \frac{\varepsilon}{2}, \\
    \frac{\beta_1\left[
    			 \frac{q_1^1}{2T_1}(\|z^*\|_2+1)^2 + \left(\tfrac{\eta_1^1}{T_1} + p_1\right)V(x_0, x^*)\right] + \sum_{k=1}^N \frac{\beta_k \sigma}{p_k c_k}}{\sum_{k=1}^{N}\beta_k} &\leq \frac{\varepsilon}{2}.
\end{align}
From Proposition \ref{pro: LO-error-proposition-stochastic}, we have that $f(x_N)-f^* \leq \varepsilon$ which implies that the resulting solution is indeed an $\varepsilon$-approximation solution. Similarly, we can also show that the consensus gap is upper-bounded by $\varepsilon$, that is $\norm{\cA x_N} \leq \varepsilon$. In the $\sigma > 0$ case, the gradient sampling complexity is:
\begin{align}
    \label{eqn: stochastic-gradient-sampling-convex}
    \sum_{k=1}^N c_k &\leq \sum_{k=1}^N \br{1+\frac{cNk}{4L^2}} = N + \frac{\sigma^2 N}{2L^2 \norm{x-x^*}^2}\sum_{k=1}^N k^2 \\
    &\leq N + \frac{\sigma N^4}{2L^2 \norm{x-x^*}^2} \leq N + \frac{800\sigma^2 \norm{x-x^*}^2}{\varepsilon^2} = O\br{N + \frac{\sigma^2}{\varepsilon^2}}
\end{align}
Notice that from proposition 2.1 in \cite{lan2021graph}, we have that $N = O\br{\frac{1}{\sqrt{\varepsilon}}}$. In addition, the number of communication rounds can be bounded as $2 \sum_{k=1}^N T_k \leq 2N + \frac{N(N+1)R\|\cA\|}{\tL} = O\br{N + \frac{1}{\varepsilon}} = O\br{\frac{1}{\varepsilon}}$.
Now, all that is left is to bound the number of LO iterations, that is $\sum_{k=1}^N \sum_{t=1}^{T_k}s_k^t$. We have:
\begin{align}
    \nonumber
    \sum_{k=1}^N \sum_{t=1}^{T_k} s_k^t &= \sum_{k=1}^N \sum_{t=1}^{T_k} \left\lfloor \frac{24(\eta_k^t+p_k)D_{\cX}^2}{\upsilon_k^t} \right\rfloor \leq \sum_{k=1}^N \sum_{t=1}^{T_k} \frac{24(\eta_k^t+p_k) D_{\cX}^2}{\varepsilon} \\
    \nonumber
    &= \sum_{k=1}^N \sum_{t=1}^{T_k} \frac{24(p_k(t+T_k)) D_{\cX}^2}{\varepsilon} \leq \sum_{k = 1}^N \frac{48p_kT_k^2D_{\cX}^2}{\varepsilon}\\
    \nonumber
    &\leq \sum_{k = 1}^N \frac{\frac{192 L}{k}\times k^2 \left\lceil \frac{R \norm{\cA}}{L} \right\rceil^2  D_{\cX}^2}{\varepsilon} \\
    &= \frac{96 L D_{\cX}^2 \left\lceil \frac{R \norm{\cA}}{L} \right\rceil^2 N(N+1)}{\varepsilon}.
\end{align}
Since we have chosen $N = O\br{\frac{1}{\sqrt{\varepsilon}}}$, the LO complexity is indeed $O\br{\frac{1}{\varepsilon^2}}$.

\textbf{Strongly-convex setting}: We proceed similarly for the strongly-convex setting. With the smoothness term of $(\eta^t_i + p_i)$ of the LO objective and diameter $D_{\cX}$ of the constraint set $\cX$, recall that the obtained error $\varepsilon_i^t$ is at most $\frac{6(\eta^t_i + p_i)L_VD_{\cX}^2}{s_i^t+1}$ where $s_i^t$ is the number of linear oracle calls at vertex $t$ and $i$-th iteration. Now, we choose:
\begin{align}
    s_i^t = \left\lfloor \frac{C}{\varepsilon} \sqrt{\frac{\beta_k}{T_k}(\eta_k^t + p_k)} \right\rfloor.
\end{align}
Since the case $k < \Delta$ is equivalent to the convex case, it is sufficient to consider $k \geq \Delta$. We have that:
\begin{align}
    \label{eqn: LO-bound-strongly-convex}
    \frac{\sum_{k=1}^N \frac{\beta_k}{T_k}(\sum_{t=1}^{T_k} \epsilon_k^t)}{\sum_{k=1}^N \beta_k} \leq \frac{\sum_{k=1}^N \frac{\beta_k}{T_k}\br{\sum_{t=1}^{T_k} \frac{6(\eta_k^t+p_k)D_{\cX}^2}{s_k^t+1}}}{\sum_{k=1}^N \beta_k} \leq \frac{\sum_{k=1}^N \frac{\beta_k}{T_k} \br{\sum_{t=1}^{T_k} \frac{6(\eta_k^t+p_k)D_{\cX}^2}{C \varepsilon \sqrt{\frac{\beta_k}{T_k}(\eta_k^t + p_k)}}}}{\sum_{k=1}^N \beta_k} 
\end{align}
for some constant $C > 0$. From $(2.14)$, we have $\eta_i^t = (p_i + \mu)(t-1) + p_iT_i, p_i = \frac{L}{1+\tau}, T_i = \left \lceil \frac{2(1+\tau)R\norm{\mathcal{A}}}{L\lambda^{\frac{i-\Delta}{2}}} \right \rceil, \beta_i = \Delta\lambda^{-(i - \Delta)}$. From this choice of parameters, notice that we have the following identities:
\begin{align}
    \sum_{k=1}^N T_k = \sum_{k=1}^N \left \lceil \frac{2(1+\tau)R\norm{\mathcal{A}}}{L\lambda^{\frac{k-\Delta}{2}}} \right \rceil \leq N + \sum_{k=1}^N \frac{2(1+\tau)R\norm{\mathcal{A}}}{L\lambda^{\frac{k-\Delta}{2}}} = N + \frac{2(1+\tau)R\norm{\mathcal{A}}}{L} \frac{\lambda^{-\frac{N-\Delta+1}{2}}-1}{\lambda^{-0.5}-1}
\end{align}
and
\begin{align}
    \label{eqn: beta-identity-strongly-convex}
    \sum_{k = \Delta }^N \beta_k = \sum_{k = \Delta }^N \Delta \lambda^{-(k-\Delta)} = \frac{\Delta(\lambda^{-(N-\Delta+1)}-1)}{\lambda^{-1}-1}.
\end{align}
Thus, combining with \eqref{eqn: LO-bound-strongly-convex}, we have:
\begin{align}
    \nonumber
    \frac{\sum_{k=1}^N \frac{\beta_k}{T_k}(\sum_{t=1}^{T_k} \epsilon_k^t)}{\sum_{k=1}^N \beta_k} &\leq  \frac{\sum_{k=1}^N \frac{\beta_k}{T_k} \br{\sum_{t=1}^{T_k} \frac{6 D_{\cX}^2 \varepsilon}{C} \sqrt{\frac{T_k (\eta_k^t + p_k)}{\beta_k}}}}{\sum_{k=1}^N \beta_k} \\
    \nonumber
    &\leq \frac{\sum_{k=1}^N \frac{\beta_k}{T_k} \br{\frac{6 D_{\cX}^2 \varepsilon}{C} \sqrt{\frac{T_k^2}{\beta_k} \sum_{t=1}^{T_k} (\eta_k^t+p_k)}}}{\sum_{k=1}^N \beta_k} \text{ from Cauchy-Schwarz} \\
    \nonumber
    &= \frac{\sum_{k=1}^N \frac{\beta_k}{T_k} \frac{6 D_{\cX}^2 \varepsilon}{C} \sqrt{\frac{T_k^2}{\beta_k} \sum_{t=1}^{T_k} ((p_k+\mu)(t-1) + p_kT_k)}}{\sum_{k=1}^N \beta_k} \\
    \nonumber
    &= \frac{\sum_{k=1}^N \frac{6 D_{\cX}^2 \varepsilon}{C} \sqrt{\beta_k \br{\frac{\br{\frac{L}{1+\tau} + \mu}T_k(T_k-1)}{2} + \frac{L\sum_{k=1}^N T_k}{1+\tau}}}}{\sum_{k=1}^N \beta_k}.
\end{align}
From here, we apply Cauchy-Schwarz once more in order to obtain the quantity $\sum_{k=1}^N \beta_k$ in the denominator. We have:
\begin{align}
    \nonumber
    &\leq \frac{\frac{6 D_{\cX}^2 \varepsilon}{C} \sqrt{\br{\sum_{k=1}^N \beta_k}\br{ \sum_{k=1}^N \frac{\br{\frac{L}{1+\tau} + \mu}(T_k^2-T_k)}{2} + \frac{N L\sum_{k=1}^N T_k}{1+\tau}}}}{\sum_{k=1}^N \beta_k} \text{ from Cauchy-Schwarz} \\
    \nonumber
    &\leq \frac{6 D_{\cX}^2 \varepsilon}{C} \sqrt{\frac{\sum_{k=1}^N \frac{\br{\frac{L}{1+\tau} + \mu}T_k^2}{2} + \frac{N L\sum_{k=1}^N T_k}{1+\tau}}{\sum_{k=1}^N \beta_k}} \\
    \nonumber
    &\leq \frac{6 D_{\cX}^2 \varepsilon}{C} \sqrt{\frac{\sum_{k=1}^N \frac{\br{\frac{L}{1+\tau} + \mu}\br{\br{\frac{2(1+\tau)R\norm{\cA}}{L\lambda^{-\frac{k-\Delta}{2}}}}^2+1}}{2} + \frac{N L\sum_{k=1}^N T_k}{1+\tau}}{\frac{\Delta(\lambda^{-(N-\Delta+1)}-1)}{\lambda^{-1}-1}}} \text{ from $\lceil x \rceil^2 \leq 4x^2+1 \forall x \geq 0$ and \eqref{eqn: beta-identity-strongly-convex}} \\
    \nonumber
    &= \frac{6 D_{\cX}^2 \varepsilon}{C} \sqrt{\frac{\sum_{k=1}^N \br{\frac{L}{1+\tau} + \mu}\frac{2(1+\tau)^2R^2\norm{\cA}^2}{L^2\lambda^{-(k-\Delta)}} + \frac{N L\sum_{k=1}^N T_k}{1+\tau} + \frac{N L}{2(1+\tau)}}{\frac{\Delta(\lambda^{-(N-\Delta+1)}-1)}{\lambda^{-1}-1}}} \\
    \nonumber
    &= \frac{6 D_{\cX}^2 \varepsilon}{C} \sqrt{\frac{\br{\frac{L}{1+\tau} + \mu}\frac{2(1+\tau)^2R^2\norm{\cA}^2}{L^2} \frac{(\lambda^{-(N-\Delta+1)}-1)}{\lambda^{-1}-1} + \frac{N L\sum_{k=1}^N T_k}{1+\tau} + \frac{N L}{2(1+\tau)}}{\frac{\Delta(\lambda^{-(N-\Delta+1)}-1)}{\lambda^{-1}-1}}} \\
    \label{eqn: LO-constant-factor-strongly-convex}
    &\leq \frac{6 D_{\cX}^2 \varepsilon}{C} \sqrt{ \br{\frac{L}{1+\tau} + \mu}\frac{2(1+\tau)^2R^2\norm{\cA}^2}{\Delta L^2} + \frac{\frac{2 N L\sum_{k=1}^N T_k}{1+\tau}}{\frac{\Delta(\lambda^{-(N-\Delta+1)}-1)}{\lambda^{-1}-1}}} \text{ since } \sum_{k=1}^N T_k \geq 1.
\end{align}
Note that we can bound the quantity in \eqref{eqn: LO-constant-factor-strongly-convex} as:
\begin{align*}
    \frac{\frac{2 N L\sum_{k=1}^N T_k}{1+\tau}}{\frac{\Delta(\lambda^{-(N-\Delta+1)}-1)}{\lambda^{-1}-1}} &= \frac{\frac{2 N L}{1 + \tau } \times \br{N + \frac{2(1+\tau)R\norm{\mathcal{A}}}{L} \frac{\lambda^{-\frac{N-\Delta+1}{2}}-1}{\lambda^{-0.5}-1}}}{\frac{\Delta(\lambda^{-(N-\Delta+1)}-1)}{\lambda^{-1}-1}} \\
    &\leq \frac{\frac{2 L N^2}{1+\tau} + \frac{4  R \norm{\cA} N \lambda^{-\frac{N-\Delta+1}{2}}}{\lambda^{-0.5}-1}}{\Delta \lambda^{-(N-\Delta)}} \\
    &\leq \frac{2LN^2}{(1+\tau)\Delta \lambda^{-(N-\Delta)}} +  \frac{4 R \norm{\cA}N}{(1-\lambda^{0.5})\Delta \lambda^{-\frac{N-\Delta}{2}}}\\
    &\leq \frac{2L}{(1+\tau)\Delta} \br{\Delta - 1 + \frac{2}{\log(\lambda^{-1})}}^2 + \frac{4 R \norm{\cA}}{(1-\lambda^{0.5})\Delta} \br{\Delta - 1 + \frac{2}{\log(\lambda^{-1})}} = C'
\end{align*}
which is upper bounded by a constant $C'$. From here, choose $C = \frac{12 D_{\cX}^2}{\sqrt{\br{\frac{L}{1+\tau} + \mu}\frac{2(1+\tau)^2R^2\norm{\cA}^2}{\Delta L^2} + C'}}$, we have that:
\begin{align}
    \label{eqn: lo-error-strongly-convex}
    \frac{\sum_{k=1}^N \frac{\beta_k}{T_k}(\sum_{t=1}^{T_k} \epsilon_k^t)}{\sum_{k=1}^N \beta_k} \leq \frac{\varepsilon}{2}.
\end{align}
In addition, also from our choice of parameter, we have the following bounds:
\begin{align}
    \frac{\beta_1\left(\frac{\eta_1^1}{T_1} + p_1\right)V(x_0, x^*) + \sum_{k=1}^N \frac{\beta_k \sigma}{p_k c_k}}{\sum_{k=1}^{N}\beta_k} &\leq \lambda^{N-\Delta}\left[ 4\tL \|x_0 - x^*\|_2^2 + \tfrac{\tL\sigma^2}{c}\right], \\
    \nonumber
    \frac{\beta_1\left[
    			 \frac{q_1^1}{2T_1}(\|z^*\|_2+1)^2 + \left(\tfrac{\eta_1^1}{T_1} + p_1\right)V(x_0, x^*)\right] + \sum_{k=1}^N \frac{\beta_k \sigma}{p_k c_k}}{\sum_{k=1}^{N}\beta_k} &\leq \lambda^{N-\Delta}\left[\tfrac{\tL}{8 R^2}(\|z^*\|_2+1)^2 \right.\\
                &\quad \left. + 4\tL \|x_0 - x^*\|_2^2 + \tfrac{\tL\sigma^2}{c}\right].
\end{align}
Thus, choose $N = \Delta + \left\lceil \log_{\lambda^{-1}} \br{\frac{10 L \norm{x_0-x^*}_2^2}{\varepsilon}} \right\rceil$, we have:
\begin{align}
    \label{eqn: LO-free-bound-primal-strongly-convex}
    \frac{\beta_1\left(\frac{\eta_1^1}{T_1} + p_1\right)V(x_0, x^*) + \sum_{k=1}^N \frac{\beta_k \sigma}{p_k c_k}}{\sum_{k=1}^{N}\beta_k} &\leq \frac{\varepsilon}{2}, \\
    \label{eqn: LO-free-bound-consensus-strongly-convex}
    \frac{\beta_1\left[
    			 \frac{q_1^1}{2T_1}(\|z^*\|_2+1)^2 + \left(\tfrac{\eta_1^1}{T_1} + p_1\right)V(x_0, x^*)\right] + \sum_{k=1}^N \frac{\beta_k \sigma}{p_k c_k}}{\sum_{k=1}^{N}\beta_k} &\leq \frac{\varepsilon}{2}.
\end{align}
Hence, combining the bounds on the primal gap in Proposition \ref{pro: LO-error-proposition-stochastic} and \eqref{eqn: lo-error-strongly-convex}, \eqref{eqn: LO-free-bound-primal-strongly-convex}, we have:
\begin{align}
    f(\overline{x}_N)-f^* &\leq \varepsilon,
\end{align}
and similarly, combining the bounds on the consensus gap in Proposition \ref{pro: LO-error-proposition-stochastic} and \eqref{eqn: lo-error-strongly-convex}, \eqref{eqn: LO-free-bound-consensus-strongly-convex}, we have:
\begin{align}
    \norm{\cA \overline{x}_N} &\leq \varepsilon,
\end{align}
which means that the obtained solution is indeed an $\varepsilon$-approximation solution. In the $\sigma > 0$ case, the gradient sampling complexity is:
\begin{align*}
    \label{eqn: stochastic-gradient-sampling-convex}
    \sum_{k=1}^N c_k &\leq \sum_{k=1}^{\Delta} \br{1+\frac{\Delta \sigma^2 k^2}{2L^2\norm{x_0-x^*}_2^2}} + \sum_{k=\Delta+1}^N \br{1 + \frac{(1+\tau)^2\Delta c \lambda^{\frac{k-\Delta}{2}}}{2L^2 \lambda^{\frac{N-\Delta}{2}}}}\\
    &\leq N + \frac{\sigma^2 \Delta^4}{2L^2 \norm{x_0-x^*}_2^2}  +  \frac{(1+\tau)^2\Delta c}{2L^2 \lambda^{\frac{N-\Delta}{2}}}\sum_{k=\Delta+1}^N \lambda^{\frac{k-\Delta}{2}} \\
    &= N + \frac{\sigma^2 \Delta^4}{2L^2 \norm{x_0-x^*}_2^2}  +  \frac{(1+\tau)^2\Delta c}{2L^2 \lambda^{\frac{N-\Delta}{2}}} \frac{\lambda^{\frac{-(N-\Delta+1)}{2}}}{\lambda^{-0.5}-1} \\
    &= N + \frac{\sigma^2 \Delta^4}{2L^2 \norm{x_0-x^*}_2^2}  +  \frac{(1+\tau)^2\Delta c}{2L^2 \lambda^{N-\Delta}(1-\lambda^{0.5})} = O\br{N + \frac{\sigma^2}{\varepsilon}}
\end{align*}
where $\Delta = O(1), \lambda^{-N} = O\br{\frac{1}{\varepsilon}}$. Now, we are left to bound the number of LO calls. We have:
\begin{align}
    \nonumber
    \sum_{k=1}^N \sum_{t=1}^{T_k} s_k^t &= \sum_{k=1}^N \sum_{t=1}^{T_k} \left\lfloor \frac{C}{\varepsilon} \sqrt{\frac{\beta_k}{T_k}(\eta_k^t + p_k)} \right\rfloor \\
    \nonumber
    &\leq \sum_{k=1}^N \sum_{t=1}^{T_k} \frac{C}{\varepsilon} \sqrt{\frac{\beta_k}{T_k}(\eta_k^t + p_k)} \\
    \nonumber
    &\leq \frac{C}{\varepsilon} \sqrt{\br{\sum_{k=1}^N T_k} \sum_{k=1}^N \frac{\beta_k}{T_k} \sum_{t=1}^{T_k} (\eta_k^t + p_k)} \text{ (C-S inequality)} \\
    \nonumber
    &= \frac{C}{\varepsilon} \sqrt{\br{\sum_{k=1}^N T_k} \sum_{k=1}^N \frac{\beta_k}{T_k} \br{\frac{\br{\frac{L}{1+\tau} + \mu}T_k(T_k-1)}{2} + \frac{L\sum_{k=1}^N T_k}{1+\tau}}} \\
    &\leq \frac{C}{\varepsilon} \sqrt{\br{\sum_{k=1}^N T_k} \frac{\br{\frac{L}{1+\tau} + \mu}\sum_{k=1}^N \beta_k T_k}{2} +  \br{\sum_{k=1}^N T_k}^2 \sum_{k=1}^N \frac{L}{1+\tau}\frac{\beta_k}{T_k}}.
\end{align}
Since we have shown above that $N = O\br{\log \frac{1}{\varepsilon}}, \sum_{k=1}^N T_k = O\br{\frac{1}{\sqrt{\varepsilon}}}, \beta_k = O\br{\lambda^{-k}}, T_k = O\br{\lambda^{-\frac{k}{2}}} \Rightarrow \sum_{k=1}^N \beta_k = O\br{\lambda^{-N}} = O\br{\frac{1}{\varepsilon}}, \sum_{k=1}^N \beta_k T_k = O\br{\lambda^{-\frac{3N}{2}}} = O\br{\frac{1}{\varepsilon^{\frac{3}{2}}}}, \sum_{k=1}^N \frac{\beta_k}{T_k} = O\br{\lambda^{-\frac{N}{2}}} = O\br{\frac{1}{\sqrt{\varepsilon}}}$. This implies that:
\begin{align}
    \sum_{k=1}^N \sum_{t=1}^{T_k} s_k^t = O\br{\frac{1}{\varepsilon^2}}
\end{align}
which means that the LO complexity is indeed $O\br{\frac{1}{\varepsilon^2}}$.
\end{proof}

\section{Ethical considerations and computing resources}
Our paper adopts a theoretical viewpoint, and the algorithms we discuss have several potential real-world applications in decentralized training. For this reason, we believe our work does not present any direct ethical and societal concerns. For our experiments, we use Python 3.9 with all algorithms implemented using \texttt{numpy} version 1.23.2. We run the experiments on a MacBook Pro with 8 GB memory. No GPU is used for the experiments.


\end{document}